\documentclass[a4paper,12pt,openany,twoside]{article}

\usepackage[utf8]{inputenc}
\usepackage[T1]{fontenc}
\usepackage{mathptmx}
\usepackage{subfig}
\usepackage{calc}
\usepackage[all]{xy}
\usepackage[centertags]{amsmath}
\usepackage{lmodern}
\usepackage{latexsym}
\usepackage{amsfonts}
\usepackage{pstricks-add}
\usepackage{amssymb}
\usepackage{amsthm}
\usepackage{multicol}
\usepackage{fancyhdr}
\usepackage{epsfig}
\usepackage{fancybox}
\usepackage{framed}
\definecolor{shadecolor}{gray}{0.9}
\usepackage{lastpage}
\usepackage{lscape}
\usepackage[francais]{varioref}
\usepackage{mathrsfs}
\usepackage{txfonts}
\usepackage{wasysym}
\usepackage{skak}
\usepackage{wallpaper}
\usepackage{eso-pic}
\usepackage{textcomp}
\usepackage[hmargin={2.8cm,2.8cm}]{geometry} 
\usepackage{lipsum}
\usepackage{bm} 

\newlength{\defbaselineskip}
\newcommand{\setlinespacing}[1]%
{\setlength{\baselineskip}{#1 \defbaselineskip}}

\newtheorem{defi}{\textbf{Definition}}[section]

\newtheorem{theorem}{\textbf{Theorem}}[section]

\newtheorem{lemma}{\textbf{Lemma}}[section]

\newtheorem{prop}{\textbf{Proposition}}[section]

\newtheorem{remark}{\textbf{Remark}}[section]

{\small\slshape\begin{flushleft}}%
	{\end{flushleft}\normalsize\upshape}

\newtheorem{cor}{\textbf{Corollary}}[section]

\newenvironment{preuv}[1][Proof.] {\noindent\textbf{#1} }{\ \rule{0.5em}{0.5em}}

\usepackage[utf8]{inputenc}
\usepackage{t1enc}
\usepackage{pstricks-add}
\usepackage{pifont}
\usepackage{ifpdf}
\ifpdf
\usepackage[pdftex=true,
hyperindex=true, colorlinks=true]{hyperref} \else
\usepackage[hypertex=true,
hyperindex=true, colorlinks=true]{hyperref} \fi

\usepackage{hyperref}
\hypersetup{
	colorlinks=true, urlcolor=black,citecolor=black,linkcolor=black,breaklinks=true}

\makeatletter

\newcommand{\Cadre}[3]{\begin{flushleft}
		\boxput*(-0.3,1){\colorbox{white}{#1}}{\setlength{\fboxsep}{8pt}
			\fbox{\begin{Bflushleft}
					\begin{minipage}{#2}
						\vspace{0.5cm}\begin{center}\par#3 \end{center} \vspace{0.5cm}
					\end{minipage}
		\end{Bflushleft}}}
	\end{flushleft}
}

\newcommand{\thechapterwords}{
	\ifcase\thechapter\or 1 \or 2\or 3 \or 4 \or
	5 \or 6 \or 7 \or 8 \or 9 \or 10 \fi}

\def\@makechapterhead#1{%
	
	\Cadre{\Huge{\scshape{\@chapapp{}\thechapterwords}}}{0.97\linewidth}{%
		\Huge \bfseries #1} \vspace{2cm} \par
	
}

\usepackage[utf8]{inputenc}
\usepackage{picinpar}
\usepackage{multicol}
\usepackage{textcomp}
\usepackage{cmap}
\usepackage{lmodern}
\usepackage{shadow}
\usepackage{makeidx}
\usepackage{mathrsfs}
\usepackage{graphics}
\usepackage{graphicx}
\usepackage{xcolor}
\usepackage{hyperref}

\begin{document}

\date{} 


\thispagestyle{empty}
\begin{center}	
\textbf{\Large{Two-scale convergence on forms in Riemannian manifolds and homogenization of an integral functional on Orlicz-Sobolev's spaces} } 
\end{center}

\vspace{1.2cm}
\begin{center}	
\footnotesize{FRANCK ARNOLD TCHINDA$^{\ddagger}$, JOEL FOTSO TACHAGO$^{\dagger}$, AND JOSEPH DONGHO$^{\dagger^{\dagger}}$}
\end{center}
\vspace{1cm}
\noindent{\textbf{Abstract :}}
We extend the concept of two-scale convergence on forms in Orlicz-Sobolev's spaces and we describe the homogenization for a family of integral functionals with convex and nonstandard growth integrands defined on the tangent bundle of a Riemannian manifold.

\vspace{0.7cm}

\noindent{\textbf{Keywords} :} 
 Two-scale convergence, forms, Riemannian manifold, Homogenization,   integral functionals, Orlicz-Sobolev's spaces.
\vspace{0.3cm}

\noindent{\textbf{2020 Mathematics Subject Classification} : 58A10, 46E30, 35B27, 49J35.}

\vspace{0.7cm}

\section{Introduction} \label{sect1} 

\par  Let $\Phi : [0,+\infty[ \rightarrow [0,+\infty[ $ be a Young function of class   $\Delta_{2}\cap \nabla_{2}$, such that its Fenchel's conjugate $\widetilde{\Phi}$ is also a Young function of class $\Delta_{2}\cap \nabla_{2}$, and $M$ an $n$-dimensional compact oriented Riemannian manifold with smooth boundary. We assume that $M$ is parallelizable.
Let $(p,\zeta_{p}) \rightarrow f(p, \zeta_{p})$ be a function from $TM$ into $\mathbb{R}$ satisfying the following properties:
\begin{itemize}
	\item[(H$_{1}$)] $f$ is measurable ;
	\item[(H$_{2}$)] $f(p,\cdot)$ is strictly convex for almost all $p \in M$ ;
	\item[(H$_{3}$)] There are two constants $c_{1}, \, c_{2} > 0$ such that 
	\begin{equation}
	c_{1} \Phi(|\zeta_{p}|) \leq f(p,\zeta_{p}) \leq c_{2}(1+ \Phi(|\zeta_{p}|))
	\end{equation}
for almost all $p \in M$ and for all $\zeta_{p} \in T_{p}M$.      
\end{itemize}
\noindent Let $U$ be an open subset of $M$ with compact closure and let $x=(x^{i})= (x^{1}, \cdots, x^{n}) : U \rightarrow \mathbb{R}^{n}$ be a coordinate system in $U$. As usual, we denote the basis on the fibers of $TU$ and $T^{\ast}U$ by $\left\{\frac{\partial}{\partial x^{i}} \right\}$ and $\left\{d x^{i} \right\}$ respectively. The coordinate system $(x^{i})$ induces the natural coordinate system $((x^{i}), (y^{i})) = (x^{1}, \cdots, x^{n}, y^{1}, \cdots, y^{n})$ on $TU$, where 
\begin{equation}
\left( x, y^{i}\frac{\partial}{\partial x^{i}} \right) \in TU \cong \left(x^{1}(p), \cdots, x^{n}(p),y^{1}, \cdots, y^{n} \right) \in \mathbb{R}^{2n}.
\end{equation}
Hence the tangent bundle $TU$ is identified with $(x^{i})(U)\times \mathbb{R}^{n}$ and the maps $f$ is denoted by $f(x,y) = f(x^{1}, \cdots, x^{n}, y^{1}, \cdots, y^{n})$. \\
Now, for each given $\epsilon > 0$, let $F_{\epsilon}$ be the integral functional defined by 
\begin{equation}\label{ras1}
F_{\epsilon}(\omega) = \int_{U} f\left( p^{\epsilon},\, ^{\sharp}d\omega_{p} \right)\, vol_{p}, \quad \omega \in H^{\Phi,d}_{0}(U),
\end{equation} 
where $p, p^{\epsilon} \in U$, $^{\sharp}d\omega_{p} \in T_{p}U$ is a vector field associate to differential 1-form $d\omega_{p}\in T^{\ast}_{p}U$,  $f$ is an integrand satisfying the above conditions $(H_{1})-(H_{3})$ and $H^{\Phi,d}_{0}(U)$ is an Orlicz-Sobolev function space on $U$ which will be specified later.  We are interested in the homogenization of the sequence of solutions of the problem
\begin{equation}\label{ras2}
\min \left\{ F_{\epsilon}(\omega) : \omega \in H^{\Phi,d}_{0}(U) \right\}.
\end{equation} 
i.e. the analysis of the asymptotic behaviour of minimizers of functionals $F_{\epsilon}$
when $\epsilon \rightarrow 0$. 
With considerations to principles of physics (see \cite{jikov}), e.g. in elasticity theory, the term $F_{\epsilon}(\omega)$ can be viewed as the energy under a deformation $\omega$ of an elastic body. 
\par In \cite{nguet1}, \textsc{NGUETSENG} develops the concept of two-scale convegence of functions followed later by \textsc{ALLAIRE} \cite{allair2}, thus allowing to homogenize the partial differential equations with periodically oscillating coefficients on a flat domain. More precisely, they establish the two-scale convergence results of functions in Lebesgue's spaces $ L^{2}(\Omega) $, where $ \Omega $ is an open bounded of $ \mathbb{R}^{n}$. 
 In \cite{tacha1}, \textsc{TACHAGO} and \textsc{NNANG} extend these results of two-scale convergence to the general case of functions defined in Orlicz's spaces 
$ L^{\Phi}(\Omega)$, where $ \Omega $ is an open bounded of $ \mathbb{R}^{n} $ and $ \Phi $ is a Young function. For application, they study the homogenization of integral functionals with convex, periodic and nonstandard growth integrands. This is the case when the functionals $F_{\epsilon}$ are of type $\int_{\Omega} f(\frac{x}{\epsilon}, Dv(x))\,dx$, where $v \in W^{1}_{0}L^{\Phi}(\Omega)$ and $W^{1}_{0}L^{\Phi}(\Omega)$ is an Orlicz-Sobolev's space of functions.  Other works on homogenization in Orcliz's spaces with flat domain can be found in \cite{tacha2, tacha3, tacha4, tacha5}. 
In \cite{tacha3}, the same authors study via the reitirated two-scale convergence method, the homogenization of integral functionals when the functionals $F_{\epsilon}$ are of type $\int_{\Omega} f(x, \frac{x}{\epsilon}, \frac{x}{\epsilon^{2}}, Dv(x))\,dx$, where $v \in W^{1}_{0}L^{\Phi}(\Omega)$.   
 In \cite{jikov}, \textsc{JIKOV V.V.} \& al. studied on a open bounded of $\mathbb{R}^{3}$, the homogenization of Maxwell's equations system in the classical formalism.  He shows that, the homogenized system obtained is well posed.
 In \cite{sango}, Sango and Woukeng propose a general method called stochastic two-scale convergence used to solve coupled periodic and stochastic homogenization problems.
 They apply it to a problem of minimzation of an integral functional where the solution depends on both variables.
 In \cite{joe}, \textsc{TACHAGO} considers the Maxwell's equations system in the classical formalism in which the solutions additionally depend on a stochastic parameter. He uses the stochastic-periodic two-scale convergence method for the homogenization of this system. Note that the work on homogenization cited so far was carried out in flat domains of $\mathbb{R}^{n}$.
  However, in order to homogenize Maxwell's equations system in their covariant formalism, \textsc{HEE. C. PAK} \cite{pak1} develops the two-scale convergence on forms called geometric two-scale convergence, which extended the two-scale convergence of \textsc{NGUETSENG} and \textsc{ALLAIRE} on the non-flat domain. More precisely, he uses the geodesics on the Riemannian manfold $M$ and establishes the two-scale convergence results in the spaces $ L^{2}(M, \Lambda^{k})$ of differential $k$-forms which are square integrable on $M$. These results will be extended later by \textsc{BACK A.} and \textsc{FRENOD E.} \cite{back2} in the spaces $ L^{r}(M, \Lambda^{k})$ of differential $k$-forms which are $r$-integrable on $M$ with $r \in [1 ; +\infty]$. Precisely, he proves the two-scale convergence compactness result (see \cite[Proposition 3]{back2}) : 
\textit{Let $M$ be an $n$-dimensional Riemannian manifold with boundary, $Y$ an $n$-dimensional Riemannian manifold without boundary. 
Assume that, $M$ is complete and $Y$ is compact and verifies the Hopf's or Mautner's conditions.   
Let $(\omega^{\epsilon})$ be a sequence in $L^{r}(M, \Lambda^{k})$ which is such that $(\omega^{\epsilon})$ is $L^{r}(M, \Lambda^{k})$-bounded and $(d^{(k)}\omega^{\epsilon})$ is $L^{r}(M, \Lambda^{k+1})$-bounded. 
Then, there exists a subsequence $(\omega^{\epsilon_{j}})$ of $(\omega^{\epsilon})$ such that :
\begin{itemize}
	\item[i)] $\omega^{\epsilon_{j}} \;\;\; \substack{2s \\ \rightharpoonup} \;\;\; \omega^{0}$ in $L^{r}(M, \Lambda^{k} L^{r}(Y))$
	\item[ii)] $d^{(k)}\omega^{\epsilon_{j}} \;\;\;  \substack{2s \\ \rightharpoonup} \;\;\; d^{(k)}\omega^{0} + d_{Y}\omega^{1}$ in $L^{r}\left(M, \Lambda^{k+1} L^{r}(Y)\right)$
\end{itemize}
with, \begin{equation}
\left\{\begin{array}{lc}
\omega^{0} \in L^{r}\left(M, \Lambda^{k} H^{1,d}(Y)\right) &  \\
\omega^{0} \in Ker(d_{Y}) &   \\
\omega^{1} \in L^{r}\left(M, \Lambda^{k+1}\, H^{1,d}(Y)\right). & 
\end{array}\right.
\end{equation}
$d^{(k)}$, $d_{Y}$ being taken in the weak sense denote the exterior and co-exterior derivative on $M$ and $Y$ respectively, $Ker(d_{Y})$ is the kernel of $d_{Y}$ and $H^{1,d}(Y)$ is the completion of $\mathcal{C}^{\infty}_{c}(Y)$ in $L^{r}(Y)$.
}
\par In this work, in order to homogenized the minimization problem (\ref{ras2}), we extend the two-scale convergence compactness results of \textsc{PAK} \cite[theorem 6.2]{pak1}, \cite[proposition 4.3]{pak1} and \textsc{BACK} \cite[theorem 4.2]{back2}, \cite[proposition 3]{back2}  to the more general spaces $ L^{\Phi}(M, \Lambda^{k})$ (see theorem \ref{math15} and proposition \ref{math5}), where $M$ is any Riemannian manifold  and $ \Phi $ is a Young function.
These spaces called Orlicz's spaces of differential forms, were apparently first examined by \textsc{YA.A. KOPYLOV} and \textsc{R.A. PANENKO} in \cite{kopy1}.
The generalization of this result to Orlicz-Sobolev's spaces (Proposition \ref{math5}) is fundamental to the proof of the main result of this work. Precisely, we show (theorem \ref{lem37}) that : \textit{ the family of minimizers of (\ref{ras2}), $(u^{\epsilon})_{\epsilon>0}$, satisfies, as $\epsilon\to 0$
	\begin{equation}
	\begin{array}{rcl}
	u^{\epsilon} & \substack{2s \\ \rightharpoonup} & u^{0} \quad \textup{in} \;\, L^{\Phi}(U)  \\
	du^{\epsilon} & \substack{2s \\ \rightharpoonup} & \textup{\textbf{d}}\textup{\textbf{u}} = du^{0} + d_{Y}u^{1} \quad \textup{in} \;\, L^{\Phi}\left(U\times Y, \Lambda^{1,0}\right) 
	\end{array} 
	\end{equation}
where $\textup{\textbf{u}}=(u^{0}, u^{1})$  is the minimizer in the space 
$H^{\Phi,d}_{0}(U)\times L^{\Phi}\left(U, \Lambda^{1}\, H^{\Phi,d}_{\#}(Y)\right)$ 
 of the global homogenized functional 
\begin{equation}
\bm{\alpha} = (\alpha^{0}, \alpha^{1}) \rightarrow \iint_{U\times Y} f(\cdot, \,^{\sharp}\textup{\textbf{d}}\bm{\alpha})\,vol_{p}vol_{q},
\end{equation}
where $^{\sharp}\textup{\textbf{d}}\bm{\alpha} = \,^{\sharp}d\alpha^{0} + \,^{\sharp}d_{Y}\alpha^{1}$.
The spaces $H^{\Phi,d}_{0}(U)$ and $L^{\Phi}\left(U, \Lambda^{1}\, H^{\Phi, d}_{\#}(Y)\right)$ will be specified later in section \ref{sect2}.
}
\par Thus, after recalling some preliminary notions on the Orlicz-Sobolev's spaces of differential forms on a Riemannian manifold in section \ref{sect2}, we establish in section \ref{sect3} the two scale-convergence results of differential forms in these spaces, and we end in section \ref{sect4} by showing the existence and uniqueness of minimizers for (\ref{ras2}) and we study the homogenization of problem for (\ref{ras2}). 

\section{Preliminaries on Orlicz-Sobolev's spaces of differential forms}\label{sect2}
\noindent{\textbf{Notations :}}
In what follows, we refer to \cite{back2}, \cite{kopy1}, \cite{iva} for notations. \\
$M$ is an $n$-dimensional oriented Riemannian manifold with smooth boundary, $\partial M$ is the boundary of $M$ and $vol_{p}$, $p\in M$ denotes the natural measure of Lebesgue on $M$. \\
$Y$ is an oriented Riemannian manifold without boundary having the same dimension as $M$.  \\
$TM$ is the tangent bundle of $M$, $T^{\ast}$ the cotangent bundle of $M$ and $|\zeta_{p}|$ denotes the norm of an element $\zeta_{p} \in T_{p}M$.\\
$\Omega^{k}(M)$, $1\leq k \leq n$ is the vector space of differential $k$-forms on  $M$. \\
$\Omega^{0}(M)$ is the space of real-valued functions on $M$.  \\
$\mathcal{C}^{\infty}(M)$ is the vector space of smooth real-valued functions on $M$. \\
$\mathcal{C}^{\infty}(M, \Lambda^{k})$ is the vector space of smooth differential $k$-forms on $M$. \\
$\mathcal{C}^{\infty}_{c}(M, \Lambda^{k})$ is the vector space of smooth differential $k$-forms with compact support in $M$. \\
$\mathcal{D'}(M, \Lambda^{k})$ is the vector space of continuous linears forms on $\mathcal{C}^{\infty}_{c}(M, \Lambda^{k})$.\\
$\Omega^{k,l}(M\times Y)$ is the vector space of differential $(k,l)$-forms on $M\times Y$. \\
$\bigwedge^{k}(T_{p}M) = \underset{k}{\underbrace{T_{p}M \wedge \cdots \wedge T_{p}M}}$ is the space of $k$-polytangent vector at $p \in M$. \\
$\bigwedge^{k}(T^{\ast}_{p}M)$, where $p \in M$ is the space of completely antisymmetric linears $k$-forms on $\bigwedge^{k}(T_{p}M)$.  \\
$\bigwedge^{k}(T^{\ast}M) = \sqcup_{p\in M} \bigwedge^{k}(T^{\ast}_{p}M)$ is the disjoint union of spaces of $k$-forms at each point $p \in M$.\\ 
$L^{r}(M, \Lambda^{k})$, where $r \in [1, +\infty]$ is the vector space of differential $k$-forms which are $r$-integrable on $M$  and $\mathcal{M}_{M}(\omega) = \frac{1}{vol(M)}\int_{M}\omega(p)vol_{p}$ is the mean value of an element $\omega \in L^{1}(M)$, whenever $M$ has a finite volume. \\
$d^{(k)} : \mathcal{C}^{\infty}(M, \Lambda^{k})\rightarrow \mathcal{C}^{\infty}(M, \Lambda^{k+1})$ is the exterior derivative operator of forms on $M$. \\
$d^{(0)}=d_{M}=d : \mathcal{C}^{\infty}(M)\rightarrow \mathcal{C}^{\infty}(M, \Lambda^{1})$ is the differential of functions on $M$. \\
$\delta^{(k)} : \mathcal{C}^{\infty}(M, \Lambda^{k})\rightarrow \mathcal{C}^{\infty}(M, \Lambda^{k-1})$ is the co-exterior derivative operator of forms on $M$ and we note  $\delta^{(0)}=\delta_{M}=\delta : \mathcal{C}^{\infty}(M, \Lambda^{1})\rightarrow \mathcal{C}^{\infty}(M).$ \\
$\star : \Omega^{k}(M)\rightarrow \Omega^{n-k}(M)$ is the Hodge star operator on $M$. \\
$\bigoplus_{i=0}^{n}\Omega^{i}(M)= \Omega(M)$ is the space of all differential forms. 

\subsection{Differential forms}
 
	A differential form of degree $k$ (or differential $k$-form) on $M$ is the application defined by :
	\begin{equation*}
	\begin{array}{rccl}
	\omega : & M & \longrightarrow & \bigwedge^{k}(T^{\ast}M) \\
	 & p & \longmapsto & \omega(p) = (p, \omega_{p})
	\end{array}
	\end{equation*} 
		We assume  that	$M$ is compact for the rest of this section. Thus any differential $k$-form on $M$ has a compact support containing in $M$. We also assume that $M$ is parallelizable i.e. for any local chart $(U,\varphi)$ of $M$ the tangent bundle $TU$ can be identified to $U\times \mathbb{R}^{n}$.
	Let $x=(x_{i})_{1\leq i\leq n}$ the coordinate system on $U$ and let $p\in U$. 
	If $\omega \in \Omega^{k}(M)$ then $\omega_{p} \in \Lambda^{k}(T^{\ast}_{p}M)$ and we have :
	\begin{equation*}
	\omega_{p} = \sum_{1\leq i_{1}< \cdots< i_{k}\leq n }\, \omega_{i_{1}, \cdots, i_{k}}(x) \; dx_{i_{1}}\wedge \cdots \wedge dx_{i_{k}}
	\end{equation*}
	where $\omega_{i_{1}, \cdots, i_{k}} : U \rightarrow \mathbb{R}$ are the functions. 
	A particular case where $\omega_{p} \in T^{\ast}_{p}M$ is given by 
	\begin{equation*}
	\omega_{p} = \omega_{1}dx_{1} + \omega_{2}dx_{2} + \cdots + \omega_{n}dx_{n}. 
	\end{equation*}
For this 1-form, one associate the vector field\, $^{\sharp}\omega$ defined by 
\begin{equation*}
^{\sharp}\omega_{p} = \omega^{1}\frac{\partial}{\partial x_{1}} + \omega^{2}\frac{\partial}{\partial x_{2}} + \cdots + \omega^{n}\frac{\partial}{\partial x_{n}},
\end{equation*}
where the coefficients of $\omega_{p}$ and $^{\sharp}\omega_{p}$ are linked by the relation 
\begin{equation*}
\omega^{i} = g^{ij}\omega_{j},
\end{equation*}
$g^{ij}$ being the elements of the inverse of metric tensor $g$ defined on $M$.  
The differential $k$-form $\omega \in \Omega^{k}(M)$ will be said to be measurable if for each point $p\in M$, its coefficients $\omega_{i_{1}, \cdots, i_{k}}$ are measurable functions on $U$ and the differential $k$-form $\omega \in \Omega^{k}(M)$ will be said to be of class $\mathcal{C}^{s}$, $s\in \overline{\mathbb{N}}$ if for each point $p\in M$, its coefficients $\omega_{i_{1}, \cdots, i_{k}}$ are the functions class $\mathcal{C}^{s}$ on $U$.  
 For any $p\in M$, we define a pointwise scalar product between two differential $k$-forms $\omega$ and $\theta$ in $\Omega^{k}(M)$ by :
\begin{equation}\label{caf1}
\langle \omega_{p},\, \theta_{p} \rangle := \int_{M} \omega_{p} \wedge \star \theta_{p}
\end{equation}
 where $\star : \Omega^{k}(M)\rightarrow \Omega^{n-k}(M)$ is the Hodge star operator and $\wedge : \Omega^{k}(M) \times \Omega^{l}(M)\rightarrow \Omega^{k+l}(M)$ the exterior product. 
 There exists an operator denoted $d^{(k)}$ from $\mathcal{C}^{\infty}(M, \Lambda^{k})$ to $\mathcal{C}^{\infty}(M, \Lambda^{k+1})$ such that : 
\begin{enumerate}
	\item[i)] $d^{(0)}=d : \mathcal{C}^{\infty}(M) \rightarrow \mathcal{C}^{\infty}(M, \Lambda^{1})$ is a differential of functions,
	\item[ii)] for all $\omega \in \mathcal{C}^{\infty}(M, \Lambda^{k})$ and $\theta \in \mathcal{C}^{\infty}(M, \Lambda^{l})$, we have $d( \omega\wedge \theta) = d\omega\,\wedge \theta + (-1)^{k}d\theta\,\wedge\omega$,
	\item[iii)] for all $\omega \in \mathcal{C}^{\infty}(M, \Lambda^{k})$,  $d\circ d \, \omega = 0$.
\end{enumerate}	
This operator is called the exterior derivative (or exterior differential) on $\mathcal{C}^{\infty}(M, \Lambda^{k})$ and if 
 $\omega_{p} = \sum_{1\leq i_{1}< \cdots< i_{k}\leq n }\, \omega_{i_{1}, \cdots, i_{k}}(x) \; dx_{i_{1}}\wedge \cdots \wedge dx_{i_{k}}$, then, 
\begin{equation*}
d\omega_{p} = \sum_{1\leq i_{1}< \cdots< i_{k}\leq n }\, \dfrac{\partial \omega_{i_{1}, \cdots, i_{k}}(x)}{\partial x_{i}}  \;dx_{i}\wedge dx_{i_{1}}\wedge \cdots \wedge dx_{i_{k}}.
\end{equation*}  
	The operator $\delta^{(k)} : \mathcal{C}^{\infty}(M, \Lambda^{k}) \rightarrow \mathcal{C}^{\infty}(M, \Lambda^{k-1})$ defined by 
	\begin{equation*}
	\delta^{(k)} = (-1)^{n(k+1)+1}\star d^{(n-k)}\star
	\end{equation*}
	is the (formal) adjoint operator of $d^{(k-1)}$ on $\bigoplus_{i=0}^{n}\mathcal{C}^{\infty}(M, \Lambda^{i})$ with respect to the scalar product $\langle \cdot \, ; \, \cdot \rangle$ defined in (\ref{caf1}).	
	$\delta^{(k)}$ is called the co-exterior derivative (or co-exterior differential) on $\mathcal{C}^{\infty}(M, \Lambda^{k})$. We will sometimes denote $\delta^{(k)}$ by $\delta^{(k)}_{M}$ and $d^{(k)}$ by $d^{(k)}_{M}$ when there is ambiguity.  
Now, we assume that $Y$ is compact and have the same dimension with $M$. A differential $(k,l)$-forms on $M\times Y$ is a map
\begin{equation*}
\begin{array}{lccr}
\omega : & M \times Y & \longrightarrow & \bigwedge^{k}(T^{\ast}M)\times \bigwedge^{l}(T^{\ast}Y) \\
& (p, q) & \longmapsto & \omega(p, q) = (p, q, \omega_{(p,q)})
\end{array}
\end{equation*} 
in which $\omega_{(p,q)}$ is a $k$-form on $M$ and a $l$-form on $Y$.
	Let  $(U, \varphi)$ be a local chart of $M$ containing the point $p$, with coordinate system $x=(x_{i})_{1\leq i\leq n}$ and $(V, \psi)$ be a local chart of $Y$ containing the point $q$, with coordinate system $y=(y_{j})_{1\leq j\leq n}$. \\
	If $\omega \in \Omega^{k,l}(M\times Y)$ then $\omega_{(p,q)} \in \Lambda^{k}(T^{\ast}_{p}M) \times \Lambda^{l}(T^{\ast}_{q}Y)$ and we have :
	\begin{equation*}
	\omega_{(p,q)} = \sum\limits_{\substack{ 1\leq i_{1}< \cdots< i_{k}\leq n \\ 1\leq j_{1}< \cdots< j_{l}\leq n }} \, \omega_{ i_{1}, \cdots, i_{k},j_{1}, \cdots, j_{l} }(x,y) \; dx_{i_{1}}\wedge \cdots \wedge dx_{i_{k}}\wedge dy_{j_{1}}\wedge \cdots \wedge dy_{j_{l}}
	\end{equation*}
	where $\omega_{i_{1}, \cdots, i_{k},j_{1}, \cdots, j_{l}} : U\times V \rightarrow \mathbb{R}$ are the functions.  
For a $\mathcal{C}^{\infty}$-mapping $f : M \rightarrow Y$ between manifolds and a differential $k$-form $\omega$ on $Y$, we define a differential $k$-form on $M$, denoted $f^{\star}(\omega)$ and call \textit{pullback} of $\omega$ by $f$ : 
\begin{equation*}
(f^{\star}(\omega))_{p}(\zeta_{1}, \cdots, \zeta_{k}) = \omega_{f(p)}(f_{\star}\zeta_{1}, \cdots, f_{\star}\zeta_{k}),
\end{equation*}
where $p \in M$, $\zeta_{1}, \cdots, \zeta_{k} \in T_{p}M$ and $f_{\star} : TM \rightarrow TY$, with $f_{\star}$ denoting the \textit{push forward}.

\subsection{Young's function}
For more results about the Young's functions, we can refer to \cite{adam}. \\
	Let $\Phi : \mathbb{R} \rightarrow  \mathbb{R}$ be a function. $\Phi$ is called a Young function if :
	\begin{itemize}
		\item[i)] $\Phi$ is continuous and convex ,
		\item [ii)] $\Phi(x) > 0$ for $x > 0$,
		\item [iii)] $\Phi(x) = 0 \; \Longleftrightarrow \; x=0$,
		\item[iv)] $\displaystyle \lim_{x \to 0} \dfrac{\Phi(x)}{x} = 0$ and $\displaystyle \lim_{x \to \infty} \dfrac{\Phi(x)}{x} = \infty$
	\end{itemize}
	The Fenchel's conjugate (or complementary function) of the Young function $\Phi$, is the Young function denoted by $\widetilde{\Phi}$ and defined by :
	\begin{equation*}
	\widetilde{\Phi}(y) = \sup_{x\geq 0} \left[ xy - \Phi(x) \right], \qquad y\geq 0.
	\end{equation*}
	A Young function $\Phi : \mathbb{R}\rightarrow \mathbb{R}$ is said to satisfy the $\Delta_{2}$-condition or $\Phi$ belongs to the class $\Delta_{2}$ for large $x$ (resp. for small $x$ or for all $x$), which is written as $\Phi \in \Delta_{2}(\infty)$ (resp. $\Phi \in \Delta_{2}(0)$ or $\Phi \in \Delta_{2}$), if there exist constants  $x_{0} > 0$, $k > 2$ such that  
	\begin{equation*}
	\Phi(2x) \leq k \Phi(x),
	\end{equation*}
	for $x \geq x_{0}$ (resp. for $0 \leq x \leq x_{0}$ or for all $x \geq 0$).	 \\
A Young function $\Phi : \mathbb{R}\rightarrow \mathbb{R}$ is said to satisfy the $\nabla_{2}$-condition or $\Phi$ belongs to the class $\nabla_{2}$ for large $x$ (resp. for small $x$ or for all $x$), which is written as $\Phi \in \nabla_{2}(\infty)$ (resp. $\Phi \in \nabla_{2}(0)$ or $\Phi \in \nabla_{2}$), if there exist constants $x_{0} > 0$, $c > 2$ such that  
	\begin{equation*}
	\Phi(x) \leq \dfrac{1}{2c} \Phi(cx),
	\end{equation*}
	for $x \geq x_{0}$ (resp. for $0 \leq x \leq x_{0}$ or for all $x \geq 0$).	 \\
	For examples, the function $x \rightarrow \frac{x^{p}}{p}$, ($p > 1$) is a Young function which satisfy $\Delta_{2}$-condition and $\nabla_{2}$-condition, its Fenchel's conjugate is the function $x \rightarrow \frac{x^{q}}{q}$, where $\frac{1}{p} + \frac{1}{q} = 1$. The function $x \rightarrow x^{p}\ln(1+x)$, ($p \geq 1$) is a Young function that satisfies $\Delta_{2}$-condition, while the Young functions $x \rightarrow x^{\ln x}$ and $x \rightarrow e^{x^{\alpha}} - 1$, ($\alpha >0$) are not of class $\Delta_{2}$.  \\
	Let a Young function $\Phi \in \Delta_{2}$. Then there are $k> 0$ and $t_{0} \geq 0$ such that 
	\begin{equation}\label{lem11}
	\widetilde{\Phi}(\phi(t)) \leq k \, \Phi(t) \quad \textup{for} \, \textup{all} \, t \geq t_{0}.
	\end{equation}	
	Let $t\rightarrow \Phi(t) = \int_{0}^{t} \phi(\tau) d\tau$ be a Young function, and let $\widetilde{\Phi}$ be the Fenchel's conjugate of $\Phi$. Then one has
	\begin{equation}\label{lem10} 
	\left\{ \begin{array}{l}
	\dfrac{t\,\phi(t)}{\Phi(t)} \geq 1 \quad (\textup{resp.} \, > 1) \; \textup{if} \, \phi \, \textup{is} \, \textup{strictly} \, \textup{increasing}  \\
	\widetilde{\Phi}(\phi(t)) \leq t\,\phi(t) \leq \Phi(2t)
	\end{array}\right.
	\end{equation}
	for all $t >0$.
	
\subsection{Orlicz-Sobolev's spaces of differential forms}
Fundamentals of Orlicz-Sobolev's spaces of differential forms can be found in \cite{kopy1}, \cite{panen1}, \cite{bour}. \\
Let $M$ and $Y$ be two $n$-dimensional Riemannian oriented compact manifold and $\Phi$ a Young function.  We describe the Orlicz spaces of differential forms. \\ 
We recall that, for any $p \in M$ in the local chart $(U, \varphi)$ with coordinate system $x=(x_{i})_{1\leq i\leq n}$, any differential $k$-form $\omega \in \Omega^{k}(M)$ can be written as 
\begin{equation*}
\omega_{p} = \sum_{1\leq i_{1}< \cdots< i_{k}\leq n }\, \omega_{i_{1}, \cdots, i_{k}}(x) \; dx_{i_{1}}\wedge \cdots \wedge dx_{i_{k}}
\end{equation*}
We denote by $vol_{p}$ the natural measure of Lebesgue associated to the Riemannian metric $g$ on $M$ and defined by $vol_{p} = \sqrt{(|g_{p}|)}dx_{1}\wedge \cdots dx_{n}$ where $|g_{p}|$ is the determinant of the Riemannian metric $g$. We  recall also that The differential $k$-form $\omega \in \Omega^{k}(M)$ is said measurable if for each point $p\in M$, its coefficients $\omega_{i_{1}, \cdots, i_{k}}$ are measurable functions on $U$.  \\  
Now, we denote :
\begin{equation*}
|\omega_{p}| := \left(\sum_{1\leq i_{1}< \cdots< i_{k} \leq n}\, (\omega_{i_{1}, \cdots, i_{k}}(x))^{2}\right)^{\frac{1}{2}}
\end{equation*}
and
\begin{equation*}
\rho_{\Phi}(\omega) := \int_{M} \Phi(|\omega_{p}|)\, vol_{p}. 
\end{equation*}
The Orlicz class of differential $k$-forms on $M$ is defined by
	\begin{equation}\label{caf5}
	 \widetilde{L}^{\Phi}(M, \Lambda^{k}) = \left\{ \omega \in \Omega^{k}(M), \,\, \textup{measurable} \; : \; \rho_{\Phi}(\omega) < +\infty  \right\}.
	\end{equation}
		 It is not a vector space. But the Orlicz space of differential $k$-forms on $M$ defined by 
	\begin{equation}
	L^{\Phi}(M, \Lambda^{k}) = \left\{ \omega \in \Omega^{k}(M), \,\, \textup{measurable} \; : \; \rho_{\Phi}(\lambda \omega) < +\infty \; for  \, some \; \lambda > 0 \right\}
		\end{equation}
 is a vector space.	 
Obviously, we have : 
\begin{equation*}
\widetilde{L}^{\Phi}(M, \Lambda^{k}) \subset L^{\Phi}(M, \Lambda^{k}).
\end{equation*}
\noindent On the manifold $M\times Y$, for any differential $(k,l)$-forms $\omega \in \Omega^{k,l}(M \times Y)$ we have :
\begin{equation*}
|\omega_{(p,q)}| := \left(\sum\limits_{\substack{1\leq i_{1}< \cdots< i_{k} \leq n \\ 1\leq j_{1}< \cdots< j_{l}\leq n }} (\omega_{ i_{1}, \cdots, i_{k},j_{1}, \cdots, j_{l}}(x,y))^{2}\right)^{\frac{1}{2}}
\end{equation*}
\begin{equation*}
\rho_{\Phi}(\omega) := \int_{M}\int_{Y} \Phi(|\omega_{(p,q)}|)\, vol_{p}vol_{q}. 
\end{equation*}
	We denote by $L^{\Phi}(M\times Y, \Lambda^{k,l})$ the Orlicz space of differential $(k,l)$-forms on $M\times Y$ defined by :
	\begin{equation*}
	L^{\Phi}(M\times Y, \Lambda^{k,l}) = \left\{ \omega \in \Omega^{k,l}(M \times Y), \,\, \textup{measurable} \; : \; \rho_{\Phi}(\lambda \omega) < +\infty \; \textup{for} \, \textup{some} \, \lambda > 0 \right\}.
	\end{equation*}
\noindent In particular, we define the subspace $L^{\Phi}\left(M, \Lambda^{k} \, \mathcal{C}^{\infty}(Y)\right)$ of $L^{\Phi}(M\times Y, \Lambda^{k,0})$ consisting of elements $\omega \in L^{\Phi}(M\times Y, \Lambda^{k,0})$ such that for a.e. $p \in M$, $\omega_{(p, \cdot)} \in \mathcal{C}^{\infty}(Y)$ and $\|\omega_{(p, \cdot)}\|_{\mathcal{C}^{\infty}(Y)} \in L^{\Phi}(M, \Lambda^{k})$ and we recall that the embbeding $L^{\Phi}\left(M, \Lambda^{k} \, \mathcal{C}^{\infty}(Y)\right) \hookrightarrow L^{\Phi}(M\times Y, \Lambda^{k,0})$ is continuous with density.
Note that	for a differential $k$-form $\omega$,
	\begin{equation*}
	\omega_{(p,q)} = \sum_{1\leq i_{1}< \cdots< i_{k}\leq n }\, \omega_{i_{1}, \cdots, i_{k}}(x,y) \; dx_{i_{1}}\wedge \cdots \wedge dx_{i_{k}} \in L^{\Phi}\left(M, \Lambda^{k} \, \mathcal{C}^{\infty}(Y)\right),
	\end{equation*}
	we can identify $\omega$ with $\bar{\omega} \in L^{\Phi}\left(M, \Lambda^{k} \, \mathcal{C}^{\infty}(Y,\Lambda^{n})\right)$
	\begin{equation*}
	\bar{\omega}_{(p,q)} = \sum_{1\leq i_{1}< \cdots< i_{k}\leq n }\, \omega_{i_{1}, \cdots, i_{k}}(x,y) \; dx_{i_{1}}\wedge \cdots \wedge dx_{i_{k}}\wedge dy_{1}\wedge \cdots \wedge dy_{n}, 
	\end{equation*}
	where $(x_{1}, \cdots, x_{n})$ and $(y_{1}, \cdots, y_{n})$ are local coordinates of $M$ and $Y$, respectively.
\\	Let now	$\Phi$ a Young function and $\widetilde{\Phi}$ its complementary.
	 $L^{\Phi}(M, \Lambda^{k})$ endowed with,
	\begin{itemize}
		\item[i)] the Luxemburg norm, 
	\begin{equation}
	\lVert \omega \rVert_{L^{\Phi}(M, \Lambda^{k})} = \inf \left\{ \lambda>0 \;\; : \;\; \rho_{\Phi}\left(\dfrac{\omega}{\lambda}\right) \leq 1 \right\}, \quad \forall \omega \in L^{\Phi}(M, \Lambda^{k}),
	\end{equation}
	\item[ii)] the Orcliz norm, 
	\begin{equation}
	\lVert \omega \rVert_{L^{(\Phi)}(M, \Lambda^{k})} = \sup \left\{ \left| \int_{M} \omega_{p} \wedge \star \theta_{p} \right| \;\; : \;\; \theta \in L^{\widetilde{\Phi}}(M, \Lambda^{k}) \; \textup{and} \; \|\theta\|_{L^{\widetilde{\Phi}}(M, \Lambda^{k})} \leq 1 \right\}, \quad \forall \omega \in L^{\Phi}(M, \Lambda^{k}),
	\end{equation}
	which makes it a Banach space.
\end{itemize}
On the other hand if we assume that $\Phi \in \Delta_{2} \cap \nabla_{2}$, Then $\mathcal{C}^{\infty}_{c}(M, \Lambda^{k})$ is dense in $L^{\Phi}(M, \Lambda^{k})$ and  we have the following embedding :
	\begin{equation*}
L^{\Phi}(M, \Lambda^{k}) \hookrightarrow L^{1}(M, \Lambda^{k}) \hookrightarrow	L^{1}_{\textup{loc}}(M, \Lambda^{k}) \hookrightarrow \mathcal{D'}(M, \Lambda^{k}).
	\end{equation*}
 Moreover the space $L^{\Phi}(M, \Lambda^{k})$ is reflexive and we have :
	\begin{equation*}
	\left( L^{\Phi}(M, \Lambda^{k}) \right)' = L^{\widetilde{\Phi}}(M, \Lambda^{k}).
	\end{equation*}
Let $p \in M$, the bilinear function, 
	\begin{equation*}
	\langle \omega_{p}, \theta_{p} \rangle_{\Phi, \widetilde{\Phi}} := \int_{M} \omega_{p} \wedge \star\theta_{p}
	\end{equation*}
	defines the association between $L^{\Phi}(M, \Lambda^{k}) \ni \omega$ and its dual $ L^{\widetilde{\Phi}}(M, \Lambda^{k}) \ni \theta$.
\\ Let $\omega \in L^{\Phi}(M, \Lambda^{k})$ and $\theta \in L^{\widetilde{\Phi}}(M, \Lambda^{k})$. For any point $p \in M$  we have Holder's inequality,
	\begin{equation*}
\left|\langle \omega_{p}, \theta_{p} \rangle_{\Phi, \widetilde{\Phi}}\right| =	\left| \int_{M} \omega_{p} \wedge \star\theta_{p} \right| \leq 2\, \lVert \omega \rVert_{L^{\Phi}(M, \Lambda^{k})} \, \lVert \theta \rVert_{L^{\widetilde{\Phi}}(M, \Lambda^{k})}.
	\end{equation*}	
	A differential $(k+1)$-form $\theta \in L^{1}_{\textup{loc}}(M, \Lambda^{k+1})$ is called a weak exterior derivative of a $k$-form $\omega \in L^{1}_{\textup{loc}}(M, \Lambda^{k})$ and we  write $d^{(k)}\omega = \theta$, if for every orientable domain $V\subset \textup{Int}M$,
	\begin{equation*}
	\int_{V} \theta_{p} \wedge u_{p} = (-1)^{k+1} \int_{V} \omega_{p} \wedge d^{(k)}\,u_{p}
	\end{equation*}
	for any $u \in \mathcal{C}^{\infty}_{c}(V, \Lambda^{n-k-1})$ and $p\in V$. 
	In the same way, we define the weak co-exterior derivative of a $k$-form $\omega \in L^{1}_{\textup{loc}}(M, \Lambda^{k})$ that we note also $\delta^{(k)}\omega$. 
\\ To end this section, we define some Orlicz-Sobolev's space of differential forms.
For $k=0 ; 1 ; \cdots ; n$,	we defined the following Orlicz-Sobolev's spaces :
	\begin{equation*}
	H^{\Phi, d^{(k)}}(M, \Lambda^{k}) = \left\{ \omega \in L^{\Phi}(M, \Lambda^{k}) \; : \; d^{(k)}\omega \in L^{\Phi}(M, \Lambda^{k+1})  \right\}
	\end{equation*}
	\begin{equation*}
	H^{\Phi, \delta^{(k)}}(M, \Lambda^{k}) = \left\{ \omega \in L^{\Phi}(M, \Lambda^{k}) \; : \; \delta^{(k)}\omega \in L^{\Phi}(M, \Lambda^{k-1})  \right\}
	\end{equation*}
	where $d^{(k)}$ and $\delta^{(k)}$ respectively denotes the exterior derivative and the co-exterior derivative on forms  in the weak sense. 
	 They have the following respective norms : 
	\begin{equation*}
	\lVert \omega \rVert_{H^{\Phi, d^{(k)}}(M, \Lambda^{k})} = \lVert \omega \rVert_{L^{\Phi}(M, \Lambda^{k})} + \lVert d^{(k)}\omega \rVert_{L^{\Phi}(M, \Lambda^{k+1})}
	\end{equation*}
	\begin{equation*}
	\lVert \omega \rVert_{H^{\Phi, \delta^{(k)}}(M, \Lambda^{k})} = \lVert \omega \rVert_{L^{\Phi}(M, \Lambda^{k})} + \lVert \delta^{(k)}\omega \rVert_{L^{\Phi}(M, \Lambda^{k-1})}
	\end{equation*}
	which makes it the Banach spaces. \\
	We denote by $H^{\Phi, d}_{\#}(M)$ the subspace of $H^{\Phi, d}(M)$ defined by
	\begin{equation*}
	H^{\Phi, d}_{\#}(M) = \left\{ \omega \in H^{\Phi, d}(M) \; : \; \mathcal{M}_{M}(\omega) = 0  \right\}
	\end{equation*}
	$H^{\Phi, d}_{\#}(M)$ is a Banach space with the following equivalent norm :
	\begin{equation*}
	\lVert \omega \rVert_{H^{\Phi, d}_{\#}(M)} = \lVert d\,\omega \rVert_{L^{\Phi}(M, \Lambda^{1})} 
	\end{equation*}

\section{Two-scale convergence of differential forms in Orlicz-Sobolev's spaces}\label{sect3} 
In order to find the homogenized functional $F$ of the sequence $F_{\epsilon}$ in (\ref{ras1}), we need to extend the notion of two-scale convergence developed in \cite{pak1} to Orlicz's spaces of differential forms $L^{\Phi}(M,\wedge^{k})$.
\subsection{Definitions}
We shall first recall the definition of two-scale convergence in Orlciz's spaces of functions on flat domains.
\begin{defi}\cite{tacha1}  \\
	Let $\Omega$ be an open bounded subset of $ \mathbb{R}^{n} $, $Y = [0 ; 1]^{n}$ the unit cube of $ \mathbb{R}^{n} $ and $\Phi$ a Young function of class $\Delta_{2}$. We note $L^{\Phi}(\Omega)$ the Orlicz space of functions on $\Omega$ and  $L^{\Phi}(\Omega \times Y_{per})$ the Orlicz space of functions on $\Omega \times Y$ which are $Y$-periodic.  \\
	A sequence of functions $(u^{\epsilon})_{\epsilon>0}$ in $L^{\Phi}(\Omega)$ is said to be weakly two-scale convergence to the function $u$ in the space $L^{\Phi}(\Omega \times Y_{per})$ if for any function $\psi \in L^{\widetilde{\Phi}}(\Omega, \mathcal{C}_{per}(Y))$, we have :  
	\begin{equation}\label{for4}
	\lim_{\epsilon\rightarrow 0} \int_{\Omega} u^{\epsilon}(x) \psi\left(x, \frac{x}{\epsilon}\right) dx = \int_{\Omega}\int_{Y} u(x,y) \psi(x,y) dx\,dy.
	\end{equation} 
\end{defi}

Now, observing formula (\ref{for4}), being given two Riemannian manifolds $M$ and $Y$, one must first define the meaning the evaluation in $\frac{x}{\epsilon}$ for differential forms on $M\times Y$. To explain this, we will use the geodesics on the manifolds $M$ and $Y$. 
\noindent We assume that $M$ is an $n$-dimensional Riemannian manifold with boundary and $Y$ is an $n$-dimensional Riemannian manifold without boundary.   
 Let $p_{0}\in M\backslash \partial M$, $q_{0}\in Y$. According to \cite{back1}, there exists an isomorphism $j$ such that $T_{p_{0}}M  \;\; \overset{j}{\cong} \;\; T_{q_{0}}Y$. 
  Moreover, if $V_{0}$ is the largest open in $T_{p_{0}}M$ containing the origin $0$, then for all $p \in \exp_{p_{0}}^{M}(V_{0})$, there exists $v \in V_{0} \subset T_{p_{0}}M$ such that
\begin{equation*}
p = \exp_{p_{0}}^{M}( v ),
\end{equation*}
where $\exp_{p_{0}}^{M}$ denotes the exponential map on $M$ in $p_{0}$.   
\noindent The two-scale convergence of differential forms follows from the application of Birkhoff's theorem\cite{bir}, Hopf's theorem\cite{hop} and the Mautner's theorem\cite{maut}. \\
Birkhoff's theorem says that for all probability space $(\Omega, \mu)$ and an ergodic flow $\varphi^{t}$, we have for $f \in L^{r}(\Omega, \mu)$,
\begin{equation*}
\lim_{T \to +\infty} \dfrac{1}{T} \int_{0}^{T} f(\varphi^{t}) dt = \int_{\Omega} f(x) d\mu.
\end{equation*}
For our concern, the geodesic flow must be ergodic on $Y$ to develop geometric two-scale convergence issues. \\
 The Hopf's theorem and the Mautner's theorem give the condtions for the geodesic flow to be ergodic. \\
  The Hopf's theorem stipulates that in a compact Riemannian manifold with a finite volume and with a negative curvative, the geodesic flow is ergodic. \\
   Mautner showed that in symmetric Riemannian manifold, the geodesic flows are also ergodic. For example, torus, hyperbolic spaces, Heinsenberg's space, the symmetric space of quaternion-K\"{a}hler are symmetric riemannian manifold.  
\noindent If the manifold $Y$ satisfies these conditions then it is geodesically complete and the
 Hopf-Rinow theorem\cite{back1} says there exists $v \in V_{0} \subset T_{p_{0}}M$ such that any point $q \in Y$ can be written as 
 \begin{equation*}
 q = \exp_{q_{0}}^{Y}(j(v)),
 \end{equation*}
 where $\exp_{q_{0}}^{Y}$ is the exponential map on $Y$ in $q_{0}$.  
\noindent  To use Birkhoff's theorem with an ergodic flow, we suppose that $M$ and $Y$ are $n$-dimensional Riemannian manifolds, moreover $M$ is assumed to be geodesically complete and possibly boundary and $Y$ assumed to be compact, without boundary and with ergodic geodesic flow i.e. verify the Mautner's condtion or Hopf's condition.  \\
 With the properties of $M$ and $Y$, for any given $p_{0}\in M$ and $q_{0} \in Y$, we define for all $p \in M$, the point $p^{\epsilon}$ which depends upon $p_{0}$, $q_{0}$ and $j$ as
\begin{equation}\label{cal1}
p^{\epsilon} = \exp_{q_{0}}^{Y}\left( \frac{1}{\epsilon}j(v) \right),
\end{equation}
for small $\epsilon > 0$ and where $v \in V_{0}$ is such that $p = \exp_{p_{0}}^{M}(v)$. \\
Having introduced this notation (\ref{cal1}), we can see that if $\psi_{(p,q)}$ is a differential $k$-form on  $M$ and a differential $0$-form on $Y$ at point $(p,q) \in M\times Y$, then $\psi_{(p,p^{\epsilon})}$ is a differential $k$-form on $M$. 
\par We have defined all the tools to do the two-scale convergence of differential forms in Orlcz's spaces.
\begin{defi} 
	Let $M$ be an $n$-dimensional Riemannian manifold with boundary, $Y$ an $n$-dimensional Riemannian manifold without boundary and $\Phi$ a Young function of class $\Delta_{2}\cap \nabla_{2}$. 
	We assume that, $M$ is complete and $Y$ is compact and verifies the Hopf's or Mautner's conditions.   
	\begin{enumerate}
\item We say that a sequence of differential forms $(\omega^{\epsilon})_{\epsilon>0}$ in $L^{\Phi}(M, \Lambda^{k})$ is weakly two-scale convergence to a form $\omega^{0}$ in $L^{\Phi}(M \times Y, \Lambda^{k,0})$, if for any form $\psi$ in $L^{\widetilde{\Phi}}(M, \Lambda^{k}\, \mathcal{C}^{\infty}(Y))$, we have   
	\begin{equation}\label{for5}
	\lim_{\epsilon\rightarrow 0} \int_{M} \omega^{\epsilon}_{p}\wedge \star \psi_{(p,p^{\epsilon})} = \int_{M}\int_{Y} \omega^{0}_{(p,q)} \wedge \star \psi_{(p,q)} vol_{q}
	\end{equation} 
	that we still write,
	\begin{equation}
\lim_{\epsilon\rightarrow 0}\;	\langle \omega_{p}^{\epsilon}, \psi_{(p,p^{\epsilon})} \rangle_{\Phi,\widetilde{\Phi}} = \langle \omega^{0}_{(p,q)}, \psi_{(p,q)}  \rangle_{\Phi,\widetilde{\Phi}}
	\end{equation}
	We will note $\omega^{\epsilon}\;\, \overset{2s}{\rightharpoonup}  \;\,\omega^{0}$ to mean it.
\item We say that a sequence of differential forms $(\omega^{\epsilon})_{\epsilon>0}$ in $L^{\Phi}(M, \Lambda^{k})$ is strongly two-scale convergence to a form $\omega^{0}$ in $L^{\Phi}(M \times Y, \Lambda^{k,0})$ if,  
\begin{equation}
\lim_{\epsilon\rightarrow 0} \lVert \omega^{\epsilon}_{p} - \omega^{0}_{(p,p^{\epsilon})} \rVert_{L^{\Phi}(M, \Lambda^{k})} = 0 
\end{equation}  	
\end{enumerate}
\end{defi}

\begin{remark}\label{lem16} 
 Assume that a sequence $(\omega_{\epsilon})$ in $ L^{\Phi}(M, \Lambda^{k})$ is weakly two-scale converge to the form $u_{0} \in   L^{\Phi}(M\times Y, \Lambda^{k,0})$. Then since $\mathcal{C}_{c}^{\infty}(M, \Lambda^{k} \mathcal{C}^{\infty}(Y))$ is dense in $L^{\widetilde{\Phi}}(M, \Lambda^{k} \mathcal{C}^{\infty}(Y))$, then (\ref{for5}) holds for every $\psi \in \mathcal{C}_{c}^{\infty}(M, \Lambda^{k} \mathcal{C}^{\infty}(Y)))$.  
\end{remark}

\subsection{Compactness theorem : case of Orlicz's spaces}
We will start by stating a lemma that we will use in the proof of the theorem.
\begin{lemma}\label{math1}
	Let $M$ an $n$-dimensional Riemannian manifold with boundary, $Y$ an $n$-dimensional Riemannian manifold without boundary and $\Phi$ a Young function of class $\Delta_{2}\cap \nabla_{2}$. 
	We assume that, $M$ is complete and $Y$ is compact and verifies the Hopf's or Mautner's conditions.   \\
	For all $ \psi \in L^{\Phi}(M, \Lambda^{k}\, \mathcal{C}^{\infty}(Y))$, We have : 
	\begin{enumerate}
		\item[i)] 
		\begin{equation} \label{jeu1}
		\lVert \psi_{(p,p^{\epsilon})} \rVert_{L^{\Phi}(M, \Lambda^{k})} \leq \lVert \psi_{(p,q)} \rVert_{L^{\Phi}(M, \Lambda^{k} \mathcal{C}^{\infty}(Y))}
		\end{equation}
		\item[ii)] 
		\begin{equation}
		\lim_{\epsilon \rightarrow 0} \lVert \psi_{(p,p^{\epsilon})} \rVert_{L^{\Phi}(M, \Lambda^{k})} = \lVert \psi_{(p,q)} \rVert_{L^{\Phi}(M \times Y, \Lambda^{k,0})}
		\end{equation}
	\end{enumerate}
\end{lemma}
\begin{preuv}
	i) Let $\epsilon > 0$, we first consider $ \psi \in \mathcal{C}^{\infty}_{c}(M, \Lambda^{k} \,\mathcal{C}^{\infty}(Y))$ then for any $p\in M$ we have :
	\begin{equation}\label{cal2}
	|\psi_{(p,p^{\epsilon})}| \leq \sup_{q\in Y} |\psi_{(p,q)}| := \lVert \psi_{p}\rVert_{\infty}
	\end{equation}
	This yields the embbeding,
	\begin{equation}
	\left\{\lambda > 0 \; : \; \int_{M} \Phi\left(\dfrac{\lVert \psi_{p}\rVert_{\infty} }{\lambda}\right)vol_{p} \leq 1  \right\} \subseteq \left\{\lambda > 0 \; : \; \int_{M} \Phi\left(\dfrac{| \psi_{(p,p^{\epsilon})}| }{\lambda}\right)vol_{p} \leq 1  \right\}
	\end{equation}
	Hence, (\ref{jeu1}) holds true for $ \psi \in \mathcal{C}^{\infty}_{c}(M, \Lambda^{k} \,\mathcal{C}^{\infty}(Y))$. The whole property follows by the density of $\mathcal{C}^{\infty}_{c}(M, \Lambda^{k} \,\mathcal{C}^{\infty}(Y))$ into $ \psi \in L^{\Phi}(M, \Lambda^{k} \,\mathcal{C}^{\infty}(Y))$. \\
	
	ii) Let $ \psi \in L^{\Phi}(M, \Lambda^{k} \,\mathcal{C}^{\infty}(Y))$. For any $\lambda >0$, we have :
	\begin{equation}
	\Phi\left(\dfrac{| \psi_{(p, q)}| }{\lambda}\right) \in L^{1}(M, \Lambda^{k}\, \mathcal{C}^{\infty}(Y)).
	\end{equation}
	It is followed  \cite[proposition 3]{back1}, that : 
	\begin{equation}
	\lim_{\epsilon\rightarrow 0} \int_{M} \Phi\left(\dfrac{| \psi_{(p, p^{\epsilon})}| }{\lambda}\right)vol_{p} = \int_{M} \int_{Y} \Phi\left(\dfrac{| \psi_{(p, q)}| }{\lambda}\right) vol_{p}vol_{q} 
	\end{equation}
	Hence, we have : 
	\begin{equation*}
	\begin{array}{rcl}
\displaystyle	\lim_{\epsilon \rightarrow 0} \lVert \psi_{(p,p^{\epsilon})} \rVert_{L^{\Phi}(M, \Lambda^{k})} & = & \displaystyle \lim_{\epsilon\rightarrow 0} \inf 	\left\{\lambda > 0 \; : \; \int_{M} \Phi\left(\dfrac{| \psi_{(p,p^{\epsilon})}| }{\lambda}\right)vol_{p} \leq 1  \right\}   \\
	& = & \inf \left\{\lambda > 0 \; : \; \displaystyle \lim_{\epsilon\rightarrow 0}\int_{M} \Phi\left(\dfrac{| \psi_{(p,p^{\epsilon})}| }{\lambda}\right)vol_{p} \leq 1  \right\}  \\
	& = & \inf \left\{\lambda > 0 \; : \; \int_{M} \int_{Y} \Phi\left(\dfrac{| \psi_{(p, q)}| }{\lambda}\right) vol_{p}vol_{q} \leq 1  \right\}  \\
	& = & \lVert \psi_{(p,q)} \rVert_{L^{\Phi}(M \times Y, \Lambda^{k,0})}.
	\end{array}
	\end{equation*}
	The proof is complete.
\end{preuv}
\begin{theorem}(compactness 1)\label{math15} \\
		Let $M$ be an $n$-dimensional Riemannian manifold with boundary, $Y$ an $n$-dimensional Riemannian manifold without boundary and $\Phi$ a Young function of class $\Delta_{2}\cap \nabla_{2}$. 
	We assume that, $M$ is complete and $Y$ is compact and verifies the Hopf's or Mautner's conditions.   \\
	Let $(\omega^{\epsilon})_{\epsilon>0}$ be a \textbf{bounded} sequence in $L^{\infty}\left([0;+\infty[ ; L^{\Phi}(M, \Lambda^{k}) \right)$. \\
	Then, there exists a subsequence $(\omega^{\epsilon_{j}})$ of $(\omega^{\epsilon})$ and a differential form $\omega^{0} \in L^{\infty}\left([0;+\infty[ ; L^{\Phi}(M\times Y, \Lambda^{k,0}) \right)$ such that for any differential form $\psi \in L^{1}\left([0;+\infty[ ; L^{\widetilde{\Phi}}(M, \Lambda^{k}\,\mathcal{C}^{\infty}(Y)) \right)$, we have : 
	\begin{equation}
	\lim_{\epsilon_{j}\rightarrow 0} \int_{0}^{+\infty} \langle \omega_{p}^{\epsilon_{j}}(t), \psi_{(p,p^{\epsilon})}(t) \rangle_{\Phi,\widetilde{\Phi}}\, dt = \int_{0}^{+\infty} \langle \omega^{0}_{(p,q)}(t), \psi_{(p,q)}(t)  \rangle_{\Phi,\widetilde{\Phi}}\, dt
	\end{equation}
\end{theorem}
\begin{preuv}
	Let $(\omega^{\epsilon}) \in L^{\infty}\left([0;+\infty[ ; L^{\Phi}(M, \Lambda^{k}) \right)$, bounded and $\psi \in L^{1}\left([0;+\infty[ ; L^{\widetilde{\Phi}}(M, \Lambda^{k} \,\mathcal{C}^{\infty}(Y)) \right)$. 
	\begin{itemize} 
	\item Consider $\mathcal{D} = L^{\widetilde{\Phi}}(M, \Lambda^{k} \,\mathcal{C}^{\infty}(Y))$ and 
	\begin{equation}
	\Upsilon_{\epsilon}(\psi) = \int_{0}^{+\infty} \langle \omega_{p}^{\epsilon_{j}}(t), \psi_{(p,p^{\epsilon})}(t) \rangle_{\Phi,\widetilde{\Phi}} dt
	\end{equation}
	\item Since $(\omega^{\epsilon})$ is bounded, then there exists a constant $c>0$ such that $$\lVert \omega^{\epsilon} \rVert_{L^{\infty}\left([0;+\infty[ ; L^{\Phi}(M, \Lambda^{k}) \right)} \leq c $$ 
	\item Since $\psi \in L^{1}\left([0;+\infty[ ; \mathcal{D} \right)$, then using the  H\"{o}lder's inequality and lemma \ref{math1}, we have : 
	\begin{equation}
	\begin{array}{rcl}
	|\Upsilon_{\epsilon}(\psi)| & \leq & \int_{0}^{+\infty} 2\,\lVert \omega_{p}^{\epsilon}(t)\rVert_{L^{\Phi}(M, \Lambda^{k})}  \lVert\psi_{(p,p^{\epsilon})}(t)\rVert_{L^{\widetilde{\Phi}}(M, \Lambda^{k})}\, dt \\
	& \leq & 2c \int_{0}^{+\infty} \lVert\psi(t)\rVert_{L^{\widetilde{\Phi}}(M; \Lambda^{k} \mathcal{C}^{\infty}(Y))} dt  \\
	& = & 2c \; \lVert \psi \rVert_{L^{1}\left([0;+\infty[ ; \mathcal{D} \right)}
	\end{array}
	\end{equation}
	\item So $\Upsilon_{\epsilon} \in \left[ L^{1}\left([0;+\infty[ ; \mathcal{D} \right) \right]'$  
	\item Moreover, since $L^{1}\left([0;+\infty[ ; \mathcal{D} \right)$ is the separable Banach space and $\Upsilon_{\epsilon}$ is bounded in $\left[ L^{1}\left([0;+\infty[ ; \mathcal{D} \right) \right]'$, then by the Banach-Alaoglu theorem \cite{bre}, there exists a subsequence $\Upsilon_{\epsilon_{j}}$ of $\Upsilon_{\epsilon}$ which weakly-$\ast$ converge to a $\Upsilon_{0}$ in $\left[ L^{1}\left([0;+\infty[ ; \mathcal{D} \right) \right]'$  
	\item Using the lemma \ref{math1} and by the dominated convergence theorem \cite{boga}, we have for any $t \in [0;+\infty[$ and for any $\psi \in L^{1}\left([0;+\infty[ ; \mathcal{D} \right)$
	\begin{equation}
	\begin{array}{rcl}
	|\Upsilon_{0}(\psi)| & = & \displaystyle \lim_{\epsilon_{j} \rightarrow 0} |\Upsilon_{\epsilon_{j}}(\psi)|  \\
	& \leq & 2c \; \int_{0}^{+\infty} \displaystyle \limsup_{\epsilon_{j} \rightarrow 0}\lVert \psi_{(p,p^{\epsilon})}(t) \rVert_{L^{\widetilde{\Phi}}(M, \Lambda^{k})} dt  \\
	& = & 2c \; \int_{0}^{+\infty} \lVert \psi(t) \rVert_{L^{\widetilde{\Phi}}(M\times Y, \Lambda^{k,0})}
	\end{array}
	\end{equation}
	\item But $L^{1}\left([0;+\infty[ ; \mathcal{D} \right)$ is dense in $L^{1}\left([0;+\infty[ ; L^{\widetilde{\Phi}}(M\times Y, \Lambda^{k,0}) \right)$, so $\Upsilon_{0}$ extends continuously into an element of 
	\begin{equation}
	\left[ L^{1}\left([0;+\infty[ ; L^{\widetilde{\Phi}}(M\times Y, \Lambda^{k,0}) \right) \right]' \equiv L^{\infty}\left([0;+\infty[ ; L^{\Phi}(M\times Y, \Lambda^{k,0}) \right)
	\end{equation}
	\item And by Riesz's theorem, we deduce that 
	\begin{equation}
	\Upsilon_{0}(\psi) = \int_{0}^{+\infty} \langle \omega^{0}(t), \psi(t) \rangle_{\Phi,\widetilde{\Phi}} dt
	\end{equation}
	with $\omega^{0} \in L^{\infty}\left([0;+\infty[ ; L^{\Phi}(M\times Y, \Lambda^{k,0}) \right)$ and $\psi \in L^{1}\left([0;+\infty[ ; L^{\widetilde{\Phi}}(M\times Y, \Lambda^{k,0}) \right)$  
	\item Thus, for any  $\psi \in L^{1}\left([0;+\infty[ ; L^{\widetilde{\Phi}}(M\times Y, \Lambda^{k,0}) \right)$ we have :
	\begin{equation}
	\begin{array}{rcl}
\displaystyle	\lim_{\epsilon_{j} \rightarrow 0} \int_{0}^{+\infty} \langle \omega_{p}^{\epsilon_{j}}(t), \psi_{(p,p^{\epsilon})}(t) \rangle_{\Phi,\widetilde{\Phi}} dt & = & \displaystyle \lim_{\epsilon_{j} \rightarrow 0} \Upsilon_{\epsilon_{j}}(\psi)  \\ 
	& = & \Upsilon_{0}(\psi)   \\
	& = & \int_{0}^{+\infty} \langle \omega^{0}(t), \psi(t) \rangle_{\Phi,\widetilde{\Phi}} dt
	\end{array}
	\end{equation}
	in particular, for all $\psi \in L^{1}\left([0;+\infty[ ; L^{\widetilde{\Phi}}(M, \Lambda^{k}\,\mathcal{C}^{\infty}(Y)) \right)$.
\end{itemize}
The proof is complete.
\end{preuv}
\begin{remark}
	Note that the compactness theorem \ref{math15} is proved in the general case with one parameter.  The simplify case stipulates that any bounded sequence $(\omega^{\epsilon})_{\epsilon> 0}$ in $ L^{\Phi}(M, \Lambda^{k}) $ admits a subsequence which  weakly two-scale converges to a $\omega^{0}$ in $L^{\Phi}(M \times Y, \Lambda^{k,0})$.
\end{remark}

\subsection{Compactness theorem : case of Orlicz-Sobolev's spaces}

\par The following lemma will be used in the proof of two-scale convergence for a derivative sequence of differential forms. The exterior derivative operator $d^{(k)}$ will be taken in the weak sense. 
\\ Let	$Y$ an $n$-dimensional Oriented compact Riemannian manifold. We reacall that, in \cite{pak1}, the Sobolev spaces $H^{1,d^{(k)}}(Y, \Lambda^{k})$ and $W^{1,d^{(k)}}(Y, \Lambda^{k})$ are respectively defined by :
\begin{equation}
H^{1,d^{(k)}}(Y, \Lambda^{k}) = \left\{ \omega \in L^{2}(Y, \Lambda^{k}) \; : \; d^{(k)}\omega \in L^{2}(Y, \Lambda^{k+1}) \right\},
\end{equation}
\begin{equation}
W^{1,d^{(k)}}(Y, \Lambda^{k}) = \left\{ \omega \in L^{1}(Y, \Lambda^{k}) \; : \; d^{(k)}\omega \in L^{1}(Y, \Lambda^{k+1}) \right\}.
\end{equation}
\begin{lemma}\label{math17}
Let	$Y$ an $n$-dimensional Oriented compact Riemannian manifold and $\Phi$ a Young function of class $\Delta_{2}\cap \nabla_{2}$.  \\
Let $\xi \in L^{\Phi}(Y, \Lambda^{1})$. Assume that for all $q \in Y$,
\begin{equation}\label{tch1}
\langle \xi_{q}, \phi_{q} \rangle_{\Phi,\widetilde{\Phi}} := \int_{Y} \xi_{q}\wedge \star\phi_{q} = 0, \quad \forall\, \phi \in \mathcal{C}^{\infty}(Y, \Lambda^{1}) \;\; \textup{such\, that} \;\; \delta^{(1)}\,\phi = 0. 
\end{equation}
Then there exists a unique $\omega \in H^{\Phi,d}_{\#}(Y)$ such that, 
\begin{equation}\label{tch2}
d\,\omega = \xi. 
\end{equation}
\end{lemma}
\begin{preuv}
	It is done in two steps.  \\
 \textbf{Step 1 :} We assume that $\xi \in \mathcal{C}^{\infty}(Y, \Lambda^{1}) \subset L^{2}(Y, \Lambda^{1})$ and verifies (\ref{tch1}).  \\
Then according to \cite[lemma 3.2]{pak1}, there exists a unique $\omega' \in H^{1,d}(Y)$ such that (\ref{tch2}) holds. \\ Consider $\omega = \omega' - \mathcal{M}_{Y}(\omega')$. Then $\omega \in H^{1,d}_{\#}(Y)$ and (\ref{tch2}) holds.  \\
Moreover, for any $p \in Y$, we have :
\begin{equation}
\int_{Y} \Phi(|d\,\omega_{p}|) \,vol_{p} = \int_{Y} \Phi(|\xi_{p}|) \,vol_{p} \leq \lVert \Phi(|\xi_{p}|) \rVert_{\infty} < +\infty. 
\end{equation}
Since $H^{1,d}_{\#}(Y) \subset W^{1,d}_{\#}(Y)$, defintion \ref{caf5} implies that, the convex set
\begin{equation}
\left\{ u \in W^{1,d}_{\#}(Y) : \int_{Y} \Phi(|d\,u_{p}|) vol_{p} < +\infty \right\},
\end{equation} 
is a subset of $H^{\Phi,d}_{\#}(Y)$. Therefore $\omega \in H^{\Phi,d}_{\#}(Y)$.  \\
 \textbf{Step 2 :} We consider the general case where $\xi \in L^{\Phi}(Y, \Lambda^{1})$ and verifies (\ref{tch1}).\\
Let $\{\varphi_{i} \}_{i\in I}$  a partition of unity subordinate to one atlas $\mathcal{A}$ for $Y$ such that it is covered by a finite number of local charts $\{(U_{i}, \kappa_{i}) \}_{i\in I}$. \\
As in \cite[page15]{iva}, we consider the \textit{pull back} of $\varphi_{i}\xi$ to $\mathbb{R}^{n}$ via $\kappa_{i}^{-1}:\mathbb{R}^{n} \rightarrow U_{i}$. This \textit{pull back} denoted $(\kappa_{i}^{-1})^{\star}(\varphi_{i}\xi)$, has  compact support and belongs to $L^{\Phi}(\mathbb{R}^{n}, \Lambda^{1})$.  \\ Thus we have that $(\kappa_{i}^{-1})^{\star}(\varphi_{i}\xi) = \sum_{1\leq j\leq n} (\kappa_{i}^{-1})^{\star}(\varphi_{i}\xi)_{j}(x)\, dx_{j}$, where $(\kappa_{i}^{-1})^{\star}(\varphi_{i}\xi)_{j} \in L^{\Phi}(\mathbb{R}^{n})$. \\
 Let $(\Psi_{m})_{m\in\mathbb{N}^{\ast}}$ be a mollifier on $\mathbb{R}^{n}$ defined by 
$\Psi_{m}(x)=m^{n}\,\Psi(mx)$, where $\Psi \in \mathcal{C}^{\infty}(\mathbb{R}^{n})$ is non negative, has a compact support in closed unit ball such that $\int \Psi(x)dx=1$. \\ For $m \in \mathbb{N}^{\ast}$, let $\omega^{i,j}_{m} = \Psi_{m} \ast (\kappa_{i}^{-1})^{\star}(\varphi_{i}\xi)_{j}$. 
Then according to \cite{adam}, we have $\omega^{i,j}_{m} \in \mathcal{C}^{\infty}(\mathbb{R}^{n})$ and \\ $\omega^{i,j}_{m} \rightarrow (\kappa_{i}^{-1})^{\star}(\varphi_{i}\xi)_{j}$ in $L^{\Phi}(\mathbb{R}^{n})$ as $m\rightarrow \infty$. Thus, $(\kappa_{i}^{-1})^{\star}(\varphi_{i}\xi)_{m} := \sum_{1\leq j\leq n} \omega^{i,j}_{m}(x) dx_{j} \in \mathcal{C}^{\infty}(\mathbb{R}^{n}, \Lambda^{1})$ and $(\kappa_{i}^{-1})^{\star}(\varphi_{i}\xi)_{m} \rightarrow (\kappa_{i}^{-1})^{\star}(\varphi_{i}\xi)$ in $L^{\Phi}(\mathbb{R}^{n}, \Lambda^{1})$ as $m\rightarrow \infty$.   \\
We now return on the manifold $Y$ by pulling back of forms $(\kappa_{i}^{-1})^{\star}(\varphi_{i}\xi)_{m}$ and $(\kappa_{i}^{-1})^{\star}(\varphi_{i}\xi)$ to $U_{i}$ via $\kappa_{i}: U_{i} \rightarrow  \mathbb{R}^{n}$. We then obtain forms $(\varphi_{i}\xi)_{m} \in \mathcal{C}^{\infty}(U_{i}, \Lambda^{1})$ and $\varphi_{i}\xi \in L^{\Phi}(U_{i}, \Lambda^{1})$ such that $(\varphi_{i}\xi)_{m} \rightarrow \varphi_{i}\xi$ in $L^{\Phi}(U_{i}, \Lambda^{1})$ as $m\rightarrow \infty$. This implies that $\sum_{i\in I} (\varphi_{i}\xi)_{m} := \xi_{m} \in \mathcal{C}^{\infty}(Y, \Lambda^{1})$, $\sum_{i\in I} \varphi_{i}\xi := \xi \in L^{\Phi}(Y, \Lambda^{1})$ and $\xi_{m} \rightarrow \xi$ in $\in L^{\Phi}(Y, \Lambda^{1})$ as $m\rightarrow \infty$. \\
Moreover, for all $\phi \in \mathcal{C}^{\infty}(Y, \Lambda^{1})$ such that $\delta^{(1)}\phi = 0$, we have $\langle \xi_{m},\,\phi \rangle = 0$. 
We therefore apply  Step1, and we conclude that there exists $\omega_{m} \in H^{\Phi,d}_{\#}(Y)$ such that
\begin{equation}\label{tch3}
 d\,\omega_{m} = \xi_{m}. 
 \end{equation} 
 In the other hand, we have $\lVert \omega_{m}-\omega_{l} \rVert_{H^{\Phi,d}_{\#}(Y)} = \lVert \xi_{m}-\xi_{l} \rVert_{L^{\Phi}(Y)}$ for all $m, l \in \mathbb{N}$ and so $(\omega_{m})$ is the Cauchy sequence in Banach space $H^{\Phi,d}_{\#}(Y)$. Therefore $\omega_{m}\rightarrow \omega$ in $H^{\Phi,d}_{\#}(Y)$ whenever $m\rightarrow \infty$. \\
 Thus, when $m\rightarrow \infty$ in (\ref{tch3}), we get that $d\,\omega = \xi$. 
\end{preuv}
\begin{prop}(compactness 2)\label{math5} \\
	Let $M$ be an $n$-dimensional Riemannian manifold with boundary, $Y$ an $n$-dimensional Riemannian manifold without boundary and $\Phi$ a Young function of class $\Delta_{2}\cap \nabla_{2}$. 
	Assume that, $M$ is complete and $Y$ is compact and verifies the Hopf's or Mautner's conditions.   \\
	Let $(\omega^{\epsilon})$ be a sequence in $L^{\Phi}(M, \Lambda^{k})$ which is such that $(\omega^{\epsilon})$ is $L^{\Phi}(M, \Lambda^{k})$-bounded and $(d^{(k)}\omega^{\epsilon})$ is $L^{\Phi}(M, \Lambda^{k+1})$-bounded. \\
	Then, there exists a subsequence $(\omega^{\epsilon_{j}})$ of $(\omega^{\epsilon})$ such that :
	\begin{itemize}
		\item[i)] $\omega^{\epsilon_{j}} \;\;\; \substack{2s \\ \rightharpoonup} \;\;\; \omega^{0}$ in $L^{\Phi}(M \times Y, \Lambda^{k,0})$
		\item[ii)] $d^{(k)}\omega^{\epsilon_{j}} \;\;\;  \substack{2s \\ \rightharpoonup} \;\;\; d^{(k)}\omega^{0} + d_{Y}\omega^{1}$ in $L^{\Phi}(M \times Y, \Lambda^{k+1,0})$
	\end{itemize}
	with, \begin{equation}
	\left\{\begin{array}{lc}
	\omega^{0} \in L^{\Phi}(M\times Y, \Lambda^{k,0}) &  \\
	\omega^{0} \in Ker(d_{Y}) &   \\
	\omega^{1} \in L^{\Phi}(M, \Lambda^{k+1}\, H^{\Phi, d}_{\#}(Y)). & 
	\end{array}\right.
	\end{equation}
	 $d^{(k)}$, $d_{Y}$ being taken in the weak sense and $Ker(d_{Y})$ is the kernel of the exterior derivative on $Y$.
\end{prop}
\begin{preuv}
Let	$(\omega^{\epsilon}) \in L^{\Phi}(M, \Lambda^{k})$ and $(d^{(k)}\omega^{\epsilon}) \in L^{\Phi}(M, \Lambda^{k+1})$. 
\begin{itemize}
	\item  Since the sequences $(\omega^{\epsilon})$ and $(d^{(k)}\omega^{\epsilon})$ are bounded, then by the previuos theorem \ref{math15}, we have :
	\begin{equation}
	\lim_{\epsilon_{j}\rightarrow 0} \int_{M} \omega^{\epsilon_{j}}_{p}\wedge \star \psi_{(p,p^{\epsilon})} = \int_{M}\int_{Y} \omega^{0}_{(p,q)} \wedge \star \psi_{(p,q)} vol_{q}
	\end{equation}
	\begin{equation}
	\lim_{\epsilon_{j}\rightarrow 0} \int_{M} d^{(k)}\omega^{\epsilon_{j}}_{p}\wedge \star \phi_{(p,p^{\epsilon})} = \int_{M}\int_{Y} \eta_{(p,q)} \wedge \star \phi_{(p,q)} vol_{q}
	\end{equation}
	with,
	\begin{equation}
	\left\{\begin{array}{lc}
	\psi \in \mathcal{C}^{\infty}_{c}(M, \Lambda^{k}\, \mathcal{C}^{\infty}(Y)) \subset L^{\widetilde{\Phi}}(M; \Lambda^{k}\, \mathcal{C}^{\infty}(Y)) &  \\
	\phi \in \mathcal{C}^{\infty}_{c}(M, \Lambda^{k+1}\, \mathcal{C}^{\infty}(Y)) \subset L^{\widetilde{\Phi}}(M, \Lambda^{k+1} \, \mathcal{C}^{\infty}(Y)) &   \\
	\omega^{0} \in L^{\Phi}(M\times Y, \Lambda^{k,0} ) &   \\
	\eta \in L^{\Phi}(M\times Y, \Lambda^{k+1,0} )  &
	\end{array}\right.
	\end{equation}
		\item  Now, since $\delta^{(k+1)}$ is the adjoint of $d^{(k)}$, we have : 
	\begin{equation}\label{rem1}
	\langle d^{(k)}\omega^{\epsilon_{j}}_{p}, \phi_{(p,p^{\epsilon})} \rangle_{\Phi,\widetilde{\Phi}} = \langle \omega^{\epsilon_{j}}_{p}, \delta^{(k+1)}_{M}\phi_{(p,p^{\epsilon})} \rangle_{\Phi,\widetilde{\Phi}} + \dfrac{1}{\epsilon} \langle \omega^{\epsilon_{j}}_{p}, \delta_{Y}\phi_{(p,p^{\epsilon})} \rangle_{\Phi,\widetilde{\Phi}}
	\end{equation}
	from where,
	\begin{equation}\label{caf2}
	\langle \omega^{\epsilon_{j}}_{p}, \delta_{Y}\phi_{(p,p^{\epsilon})} \rangle_{\Phi,\widetilde{\Phi}} = \epsilon \left( \langle d^{(k)}\omega^{\epsilon_{j}}_{p}, \phi_{(p,p^{\epsilon})} \rangle_{\Phi,\widetilde{\Phi}} - \langle \omega^{\epsilon_{j}}_{p}, \delta^{(k+1)}_{M}\phi_{(p,p^{\epsilon})} \rangle_{\Phi,\widetilde{\Phi}} \right).
	\end{equation}
		\item  Thus, passing to the two-scale limits in the equation (\ref{caf2}) , we have for all $\phi \in \mathcal{C}^{\infty}_{c}(M, \Lambda^{k+1}\, \mathcal{C}^{\infty}(Y))$ :
	\begin{equation}
	\begin{array}{rrll}	
\langle \omega^{0}, \delta_{Y}\phi \rangle_{\Phi,\widetilde{\Phi}} = 0 	& \Rightarrow & \langle d_{Y}\omega^{0}, \phi \rangle_{\Phi,\widetilde{\Phi}} = 0,	 &   \\
 &	\Rightarrow & d_{Y}\omega^{0} = 0  &  \\
 &	\Rightarrow & \omega^{0} \in Ker(d_{Y}).
	\end{array}
	\end{equation}
	\item  Moreover, let  $ \phi \in \mathcal{C}^{\infty}_{c}(M, \Lambda^{k+1}\, \mathcal{C}^{\infty}(Y))$ such that $\delta_{Y}\phi = 0$, the equation (\ref{rem1}) becomes :
	\begin{equation}\label{caf3}
\langle \omega^{\epsilon_{j}}_{p}, \delta^{(k+1)}_{M}\phi_{(p,p^{\epsilon})} \rangle_{\Phi,\widetilde{\Phi}} =	\langle d^{(k)}\omega^{\epsilon_{j}}_{p}, \phi_{(p,p^{\epsilon})} \rangle_{\Phi,\widetilde{\Phi}} 
	\end{equation}
	then passing to the two-scale limits in (\ref{caf3}) and knowing that $\omega^{0}$ does not depend on $Y$, we obtain : 
	\begin{equation}\label{caf4}
	\begin{array}{rll}
\langle \omega^{0}, \delta^{(k+1)}_{M}\phi \rangle_{\Phi,\widetilde{\Phi}} = \langle \eta, \phi \rangle_{\Phi,\widetilde{\Phi}}  &	\Rightarrow & \langle d^{(k)}\omega^{0}, \phi \rangle_{\Phi,\widetilde{\Phi}} = \langle \eta, \phi \rangle_{\Phi,\widetilde{\Phi}}	   \\
	&\Rightarrow & \langle \eta - d^{(k)}\omega^{0}, \phi \rangle_{\Phi,\widetilde{\Phi}} = 0
	\end{array}
	\end{equation}
		\item  By the lemma \ref{math17}, since (\ref{caf4}) holds for all $ \phi \in \mathcal{C}^{\infty}_{c}(M, \Lambda^{k+1}\, \mathcal{C}^{\infty}(Y)) $ such that $\delta_{Y}\phi = 0$, we conclude that there exists $\omega^{1} \in L^{\Phi}(M, \Lambda^{k+1}\, H^{\Phi, d}_{\#}(Y))$ such that 
	\begin{equation}
	\eta - d^{(k)}\omega^{0} = d_{Y}\omega^{1} \quad i.e. \quad \eta = d^{(k)}\omega^{0} + d_{Y}\omega^{1}.	
	\end{equation}
\end{itemize}
The proof is complete.
\end{preuv}  

\section{Application to the homogenization of an integral functional with convex integrand and Nonstandard growth}\label{sect4}

\par In this section, $M$ denotes an $n$-dimensional compact Riemannian manifold with smooth boundary and $Y$ an $n$-dimensional Riemannian manifold without boundary and $\Phi$ a Young function of class $\Delta_{2}\cap \nabla_{2}$. 
We assume that, $M$ is complete and parallelizable, and $Y$ is compact and verifies the Hopf's or Mautner's conditions.
Let $\Phi : [0,+\infty[ \rightarrow [0,+\infty[ $ be a Young function of class   $\Delta_{2}\cap \nabla_{2}$, such that its Fenchel's conjugate $\widetilde{\Phi}$ is also a Young function of class $\Delta_{2}\cap \nabla_{2}$.\\
Let $(p,\zeta_{p}) \rightarrow f(p, \zeta_{p})$ be a function from $TM$ into $\mathbb{R}$ satisfying the following properties:
\begin{itemize}
	\item[(H$_{1}$)] $f$ is measurable ;
	\item[(H$_{2}$)] $f(p,\cdot)$ is strictly convex for almost all $p \in M$ ;
	\item[(H$_{3}$)] There are two constants $c_{1}, \, c_{2} > 0$ such that 
	\begin{equation}
	c_{1} \Phi(|\zeta_{p}|) \leq f(p,\zeta_{p}) \leq c_{2}(1+ \Phi(|\zeta_{p}|))
	\end{equation}
	for almost all $p \in M$ and for all $\zeta_{p} \in T_{p}M$.      
\end{itemize}
Let $U$ be an open subset of $M$ with compact closure and let $x=(x^{i})= (x^{1}, \cdots, x^{n}) : U \rightarrow \mathbb{R}^{n}$ be a coordinate system in $U$.
For each given $\epsilon > 0$, let $F_{\epsilon}$ be the integral functional defined by 
\begin{equation}\label{lem13}
F_{\epsilon}(\omega) = \int_{U} f\left( p^{\epsilon}, \, ^{\sharp}d\omega_{p} \right)\, vol_{p} \quad \omega \in H^{\Phi,d}_{0}(U),
\end{equation} 
where $p, p^{\epsilon} \in U$,\, $^{\sharp}d\omega_{p} \in T_{p}U$ is a vector field associated to differential 1-form $d\omega_{p}\in T^{\ast}_{p}U$,  $f$ is an integrand satisfying above conditions $(H_{1})-(H_{3})$ and $H^{\Phi,d}_{0}(U)$ is the Orlicz-Sobolev function space on $U$.  We are interested in the homogenization of the sequence of solutions of the problem
\begin{equation}\label{lem26}
\min \left\{ F_{\epsilon}(\omega) : \omega \in H^{\Phi,d}_{0}(U) \right\}.
\end{equation} 
i.e. the analysis of the asymptotic behaviour of minimizers of functionals $F_{\epsilon}$
when $\epsilon \rightarrow 0$. 

\subsection{Setting of the problem : Existence and uniqueness of minimizers}

Let us begin by the following lemma.
\begin{lemma}\label{lem9}
	Let $(p,\alpha_{p}) \rightarrow g(p,\zeta_{p})$ be a function from $TM$ into $\mathbb{R}$ satisfying the following properties :
	\begin{itemize}
		\item[(i)] $g$ is measurable ;
		\item[(ii)] $g(p,\cdot)$ is strictly convex for almost all $p \in M$ ; 
		\item[(iii)] There exists  $c_{0} > 0$ such that 
		\begin{equation} \label{lem8}
		| g(p,\zeta_{p})| \leq c_{0}(1+ \Phi(|\alpha_{p}|))
		\end{equation}
		for almost all $p \in M$ and for all $\zeta_{p} \in T_{p}M$.
	\end{itemize}
	Then $g$ is continuous in the second argument. Precisely one has
	\begin{equation}\label{lem7}
	| g(p,\zeta_{p}) - g(p,\theta_{p}) | \leq c'_{0} \dfrac{1 + \Phi(2(1+ |\zeta_{p}| + |\theta_{p}|))}{1+ |\zeta_{p}| + |\theta_{p}|} |\zeta_{p} - \theta_{p} |
	\end{equation}
	for almost all $p \in M$ and for all $\zeta_{p},\,\theta_{p} \in T_{p}M$, with $c'_{0} = 2Nc_{0}$.
\end{lemma}
\begin{preuv}
	Let $(U,\phi)$ a local chart of $M$ containing $p$ and $(x_{1}, \cdots, x_{n})=x$, the local coordinate system on $U$. 
	Since $M$ is parallelizable, we can identify $TU \cong U\times \mathbb{R}^{n}$.
	Hence, we can write $\zeta_{p} = (\zeta^{i})= \zeta$ and $\theta_{p} = (\theta^{i})=\theta$. We will prove that 
	 \begin{equation}\label{lem70}
	 | g(x,\zeta) - g(x,\theta) | \leq c'_{0} \dfrac{1 + \Phi(2(1+ |\zeta| + |\theta|))}{1+ |\zeta| + |\theta|} |\zeta - \theta|.
	 \end{equation}
	To achieve this, let us define the vector $\mu^{k}$,  $0 \leq k \leq n$, by 
	\begin{equation*}
	\mu^{0} = \zeta,\, \mu^{n}=\theta, \;\, \textup{and} \;\, \mu^{k}= (\zeta^{1}, \cdots, \zeta^{k}, \theta^{k+1},\cdots \theta^{n}) \;\; \textup{for} \;\, 1\leq k \leq n-1.
	\end{equation*}
	Now, observing that for all $x\in \mathbb{R}^{n}$ the equation of the straight line passing by both points $(\mu^{k}, g(x,\mu^{k}))$ and $(\mu^{k+1}, g(x,\mu^{k+1}))$ is 
	\begin{equation}
	z = (\gamma - \zeta_{k+1}) \dfrac{g(x, \mu^{k+1}) - g(x, \mu^{k})}{\theta_{k+1} - \zeta_{k+1}} + g(x, \mu^{k}), \quad \gamma \in \mathbb{R},
	\end{equation}
	then putting  $\widehat{\gamma} = (\mu^{1},\cdots, \mu^{k},\zeta^{k+2},\cdots, \zeta^{n})$, we have, thanks to the convexity, 
	\begin{equation}
	g(x, \widehat{\gamma}) \geq (\gamma - \zeta_{k+1}) \dfrac{g(x, \mu^{k+1}) - g(x, \mu^{k})}{ \theta_{k+1} - \zeta_{k+1}} + g(x, \mu^{k}), 
	\end{equation}
	for all $\gamma$ outside  the interval $(\theta_{k+1} , \zeta_{k+1})$ if $\zeta_{k+1} > \theta_{k+1}$ or $(\zeta_{k+1},\theta_{k+1})$ if $\theta_{k+1} > \zeta_{k+1}$. \\
	Taking $\gamma = t\theta_{k+1} + (1-t)\zeta_{k+1}$ with $t \geq 1$, so that $\widehat{\gamma} = t\mu^{k+1} + (1-t)\mu^{k}$, we get 
	\begin{equation}
	g(x,t\mu^{k+1} + (1-t)\mu^{k}) \geq t [g(x, \mu^{k+1}) - g(x, \mu^{k})] + g(x, \mu^{k}),
	\end{equation}
	that is, for all $t \geq 1$, 
	\begin{equation}\label{lem6}
	g(x, \mu^{k+1}) - g(x, \mu^{k}) \leq \dfrac{g(x,t\mu^{k+1} + (1-t)\mu^{k}) - g(x, \mu^{k})}{t}.
	\end{equation}
	But it is easy to see that $|\mu^{k}| \leq 1 + |\zeta| + |\theta|$ and $|\mu^{k}+ t(\mu^{k+1} -\mu^{k})| \leq 1 + |\zeta| + |\theta| + t|\zeta - \theta|$, then taking $t_{0} = \dfrac{1 + |\zeta| + |\theta| }{|\zeta - \theta|} \geq 1$ yields $|\mu^{k}+ t(\mu^{k+1} -\mu^{k})| \leq 2 (1 + |\zeta| + |\theta|)$. This being so, according to hypotheses (iii) we have  
	\begin{equation}
	\dfrac{g(x,t\mu^{k+1} + (1-t)\mu^{k}) - g(x, \mu^{k})}{t_{0}} \leq 2 c_{0} \dfrac{1 + \Phi(2(1+ |\zeta| + |\theta|))}{1+ |\zeta| + |\theta|} |\zeta - \theta|.
	\end{equation}
	Then summing each member of (\ref{lem6}) on $k = 0, \cdots, n-1$, and reversing the role of $\zeta$ and $\theta$ it follows (\ref{lem7}).
\end{preuv}

Let  $\alpha = (\alpha^{i}) \in \mathcal{C}_{c}^{\infty}(\overline{U})^{n}$. By (\ref{lem8}), it is an easy task to check that the function $(p, q) \rightarrow g(q, \alpha_{p})$ of $U\times Y$ into $\mathbb{R}$ belongs to $\mathcal{C}_{c}^{\infty}(\overline{U}; L^{\infty}(Y))$. Hence for each $\epsilon > 0$, the function $p \rightarrow g(p^{\epsilon}, \alpha_{p})$ of $U$ into $\mathbb{R}$, denoted by $g^{\epsilon}(\cdot, \alpha)$ is well defined as an element of $L^{\infty}(U)$, and we have the following proposition and corollary.

\begin{prop}\label{lem14}
	Given $\epsilon >0$, the transformation $\alpha \rightarrow g^{\epsilon}(\cdot, \alpha)$ of $\mathcal{C}_{c}^{\infty}(\overline{U})^{n}$ into $L^{\infty}(U)$ extends by continuity to a mapping, still denoted by $\alpha \rightarrow g^{\epsilon}(\cdot, \alpha)$ of
	$L^{\Phi}(U)^{n}$ into $L^{1}(U)$ with the property 
	\begin{equation}\label{lem12}
	\parallel g^{\epsilon}(\cdot, \alpha) - g^{\epsilon}(\cdot, \beta) \parallel_{L^{1}(U)} \leq c\, \left( \|1\|_{\widetilde{\Phi},U} + \parallel \phi( 1+|\alpha|+ |\beta|)\parallel_{\widetilde{\Phi},U} \right)\parallel \alpha-\beta  \parallel_{L^{\Phi}(U)^{n}}
	\end{equation}
	for all $\alpha,\,\beta \in L^{\Phi}(U)^{n}$.
\end{prop}
\begin{preuv}
Let $p\in \overline{U}$, $q\in Y$ and $(x_{1}, \cdots, x_{n})=x$ the locale coordinate system on $U$. Then for $\alpha,\,\beta \in \mathcal{C}_{c}^{\infty}(\overline{U})^{n}$, we have $\alpha_{p}= (\alpha^{i}(x)) \in \mathbb{R}^{n}$, $\beta_{p}= (\beta^{i}(x)) \in \mathbb{R}^{n}$. For any fixed $\epsilon >0$ Putting 
	\begin{equation}
	\psi(p,q) = c'_{0} \dfrac{1 + \Phi(2(1+ |\alpha_{p}| + |\beta_{p}|))}{1+ |\alpha_{p}| + |\beta_{p}|} |\alpha_{p} - \beta_{p} | - |g(q, \alpha_{p}) - g(q, \beta_{p})|,
	\end{equation}
	then $\psi \in \mathcal{C}^{\infty}(\overline{U}; L^{\infty}(Y))$ with $\psi(p,q) \geq 0$ for all $p \in \overline{U}$ and for almost all $q \in Y$ (see lemma \ref{lem9}), so that $\psi(p,p^{\epsilon}) \geq 0$, that is 
	\begin{equation}
	|g^{\epsilon}(\cdot, \alpha) - g^{\epsilon}(\cdot, \beta)| \leq c'_{0} \dfrac{1 + \Phi(2(1+ |\alpha| + |\beta|))}{1+ |\alpha| + |\beta|} |\alpha - \beta| \quad \textup{a.e. \, in} \, U. 
	\end{equation}
	Since $\Phi \in \Delta_{2}$, there are $k, \delta >0$ such that $\Phi(2t) \leq k \Phi(t)$ for all $t\geq \delta$. Letting $U_{\delta} = \{ p \in U \, : \, 1+ |\alpha_{p}| + |\beta_{p}| \geq \frac{\delta}{2}  \}$ and $U^{\delta} = U\backslash U_{\delta}$, one has 
	\begin{equation}
	|g^{\epsilon}(\cdot, \alpha) - g^{\epsilon}(\cdot, \beta)| \leq c'_{0}(1 + \Phi(\delta))|\alpha - \beta|   \quad \textup{a.e. \, in} \, U^{\delta}
	\end{equation}
	and, appeals to the first inequality in (\ref{lem10}),
	\begin{equation}
	|g^{\epsilon}(\cdot, \alpha) - g^{\epsilon}(\cdot, \beta)| \leq c'_{0}(1 + k\phi(1+ |\alpha| + |\beta| ))|\alpha - \beta|   \quad \textup{a.e. \, in} \, U_{\delta},
	\end{equation}
	where $\phi$ is the right derivative of $\Phi$. Hence we get 
	\begin{equation}
	|g^{\epsilon}(\cdot, \alpha) - g^{\epsilon}(\cdot, \beta)| \leq c''_{0}(1 + \phi(1+ |\alpha| + |\beta| ))|\alpha - \beta|   \quad \textup{a.e. \, in} \, U,
	\end{equation}
	where $c''_{0} = c'_{0}\max \{ 1+ \Phi(\delta), k \}$. 
	Integrating on $U$, and applying the generalized H\"{o}lder's inequality, we get results for $\alpha,\,\beta \in \mathcal{C}_{c}^{\infty}(U)^{n}$. 
	Therefore, since $\mathcal{C}_{c}^{\infty}(\overline{U})^{n}$ is dense in $L^{\Phi}(U)^{n}$, the proposition follows by extension by continuity. 
\end{preuv}

In proposition \ref{lem14}, if we take $g=f$, then we have the following corallary. 

\begin{cor}\label{lem22}
	Under the hypothesis $(H_{1})-(H_{3})$, given $\omega \in H^{\Phi,d}(U)$, the function 
	$p \rightarrow f(p^{\epsilon}, \, ^{\sharp}d\omega_{p})$ of  $U$ into $\mathbb{R}$, denoted $f^{\epsilon}(\cdot, \, ^{\sharp}d\omega)$ is well defined as an element of $L^{1}(U)$. Moreover we have 
	\begin{equation}\label{lem15}
	c_{1} \parallel\,^{\sharp}d\omega \parallel_{L^{\Phi}(U)^{n}} \leq \parallel f^{\epsilon}(\cdot, \, ^{\sharp}d\omega) \parallel_{L^{1}(U)} \leq c'_{2}(1 + \parallel\, ^{\sharp}d\omega \parallel_{L^{\Phi}(U)^{n}}),
	\end{equation}	
	for all $\omega \in H^{\Phi,d}(U)$, where $c'_{2} = c_{2}\max(1, |U|)$ with $|U| = \int_{U} vol_{p}$. 
\end{cor}

We are now able to prove the existence of a minimizer on $H^{\Phi,d}_{0}(U)$, for each $\epsilon > 0$ of integral functional $F_{\epsilon}$ (see (\ref{lem13})). 
\begin{theorem}
	For $\epsilon > 0$, there exists a unique $u^{\epsilon} \in H^{\Phi,d}_{0}(U)$ that minimizes $	F_{\epsilon}$ on  $H^{\Phi,d}_{0}(U)$, i.e.,
	\begin{equation}
	F_{\epsilon}(u^{\epsilon}) = \min \left\{ F_{\epsilon}(\omega) : \omega \in H^{\Phi,d}_{0}(U) \right\}.
	\end{equation}
\end{theorem}
\begin{preuv}
	Let $\epsilon > 0$ be fixed. Thanks to proposition \ref{lem14} (with $g=f$), there is a constant $c > 0$ such that 
	\begin{equation}
	\left| F_{\epsilon}(\omega) - F_{\epsilon}(\upsilon)  \right| \leq c\, \left( \|1\|_{\widetilde{\Phi},U} + \parallel \phi( 1+|^{\sharp}d\omega|+ |^{\sharp}d\upsilon|)\parallel_{\widetilde{\Phi},U} \right)\parallel \,^{\sharp}d\omega - \,^{\sharp}d\upsilon\parallel_{L^{\Phi}(U)^{n}},
	\end{equation} 
	for all $\omega, \upsilon \in H^{\Phi,d}_{0}(U)$, so that $F_{\epsilon}$ is continuous. With this in mind, since $F_{\epsilon}$ is strictly convex (see $(H_{2})$), and coercive (see the left-hand side inequality in (\ref{lem15})), there exists a unique $u^{\epsilon}$ that minimizes $F_{\epsilon}$ on $H^{\Phi,d}_{0}(U)$.
\end{preuv} 

In the following we investigate to the limiting behaviour, as $\epsilon\to 0$, of the sequence of minimizers  of $u^{\epsilon}$.

\subsection{Preliminary results} 
Let $\alpha=(\alpha^{i})  \in  \mathcal{C}^{\infty}(Y)^{n}$.
The function $q \rightarrow f(q, \alpha_{q})$ of $Y$ into $\mathbb{R}$ denoted by $f(\cdot, \alpha)$ belongs to $L^{\infty}(Y)$. This being, using lemma \ref{lem9}, for every $\alpha \in \mathcal{C}_{c}^{\infty}(\overline{U}, \mathcal{C}^{\infty}(Y))$  the function
$p \rightarrow f(\cdot, \alpha_{(p,\cdot)})$ from $\overline{U}$ to $L^{\infty}(Y)$ still denoted by $f(\cdot, \alpha)$ lies in  $\mathcal{C}_{c}^{\infty}(\overline{U}; L^{\infty}(Y))$, with
\begin{equation}\label{lem19}
|f(\cdot, \alpha) - f(\cdot, \beta)| \leq c'_{2} \left(  1 + \dfrac{\Phi(2(1+ |\alpha| + |\beta|))}{1+ |\alpha| + |\beta|}\right) |\alpha - \beta| \quad \textup{a.e. \, in} \,\; \overline{U}\times Y 
\end{equation}
for all $\alpha, \beta  \in \mathcal{C}_{c}^{\infty}(\overline{U}, \mathcal{C}^{\infty}(Y))$. Therefore, for fixed $\epsilon >0$, one defines the function $p \rightarrow f(p^{\epsilon}, \alpha_{(p,p^{\epsilon})})$ of $U$ into $\mathbb{R}$, denoted by $f^{\epsilon}(\cdot, \alpha^{\epsilon})$  as an element of $L^{\infty}(U)$.
\begin{prop}\label{lem23}
	Suppose that $(H_{1})-(H_{3})$ hold. For every $\alpha  \in \mathcal{C}_{c}^{\infty}(\overline{U}, \mathcal{C}^{\infty}(Y))^{n}$ one has 
	\begin{equation}\label{lem17}
	\lim_{\epsilon \to 0} \int_{U} f(p^{\epsilon}, \alpha_{(p,p^{\epsilon})})  vol_{p} = \iint_{U\times Y} f(q, \alpha_{(p,q)}) vol_{p}vol_{q}.
	\end{equation}
	Futhermore, the mapping  $\alpha \rightarrow f(\cdot, \alpha)$ of $\mathcal{C}_{c}^{\infty}(\overline{U}, \mathcal{C}^{\infty}(Y))^{n}$ into $L^{1}(U\times Y)$ extends by continuity to a mapping, still denoted by $\alpha \rightarrow f(\cdot, \alpha)$, of
	$L^{\Phi}(U\times Y)^{n}$ into $L^{1}(U\times Y)$ with the property 
	\begin{equation}
	\parallel f(\cdot, \alpha) - f(\cdot, \beta) \parallel_{L^{1}(U\times Y)} \leq c\, \left( \|1\|_{\widetilde{\Phi},U\times Y} + \parallel \phi( 1+|\alpha|+ |\beta|)\parallel_{\widetilde{\Phi},U\times Y} \right)\parallel \alpha -\beta  \parallel_{L^{\Phi}(U\times Y)^{n}}
	\end{equation}
	for all $\alpha, \beta  \in L^{\Phi}(U\times Y)^{n}$.
\end{prop}

\begin{preuv}
	Let	$\alpha  \in \mathcal{C}_{c}^{\infty}(\overline{U}, \mathcal{C}^{\infty}(Y))$. Since $f(\cdot, \alpha)$ belongs to $\mathcal{C}_{c}^{\infty}(\overline{U}, L^{\infty}(Y))$ and the sequence of functions $(\omega^{\epsilon})$ defined by $\omega^{\epsilon}_{p} =1$, a.e. $p\in \overline{U}$ is weakly two-scale convergent in $L^{\Phi}(U)$ to 1, the convergence result (\ref{lem17}) is a consequence of remark \ref{lem16}. 
	On the other hand, since $\Phi \in \Delta_{2}$, arguing as in the proof of proposition \ref{lem14}, there is a constant $c= c(c_{2},Q,\Omega,\Phi)$ such that  
	\begin{equation}
		\parallel f(\cdot, \alpha) - f(\cdot, \beta) \parallel_{L^{1}(U\times Y)} \leq c\, \left( \|1\|_{\widetilde{\Phi},U\times Y} + \parallel \phi( 1+|\alpha|+ |\beta|)\parallel_{\widetilde{\Phi},U\times Y} \right)\parallel \alpha -\beta  \parallel_{L^{\Phi}(U\times Y)^{n}}
	\end{equation}
for all $\alpha, \beta  \in \mathcal{C}_{c}^{\infty}(\overline{U}, \mathcal{C}^{\infty}(Y))^{n}$. We end the proof by using argument of continuity and density of $\mathcal{C}_{c}^{\infty}(\overline{U}, \mathcal{C}^{\infty}(Y))^{n}$ in $L^{\Phi}(U, \mathcal{C}^{\infty}(Y))^{n}$.
\end{preuv}
\begin{cor}\label{lem24}
	Let 
	\begin{equation}
	\omega^{\epsilon}_{p} = \omega^{0}_{p} + \epsilon \omega^{1}_{(p,p_{\epsilon})},
	\end{equation}
	with $\epsilon>0$, $p \in U$, $\omega^{0} \in \mathcal{C}_{c}^{\infty}(U)$ and $\omega^{1} \in \mathcal{C}_{c}^{\infty}(U)\otimes \mathcal{C}^{\infty}(Y)$ . Then 
	\begin{equation}
	\lim_{\epsilon \to 0} \int_{U} f(p^{\epsilon}, \, ^{\sharp}d\omega^{\epsilon}_{p}) vol_{p} = \iint_{U\times Y} f(p, \, ^{\sharp}d\omega^{0}_{p} + \, ^{\sharp}d_{Y}\omega^{1}_{(p,q)}) vol_{p}vol_{q}.
	\end{equation}
\end{cor}
\begin{preuv}
	As $d\omega^{1}_{(p,p_{\epsilon})} = d\omega^{1}_{(p,p_{\epsilon})} + \frac{1}{\epsilon} d_{Y}\omega^{1}_{(p,p_{\epsilon})}$, we have  $d\omega^{\epsilon}_{p} = d\omega^{0}_{p} + \epsilon d\omega^{1}_{(p,p_{\epsilon})} + d_{Y}\omega^{1}_{(p,p_{\epsilon})}$. \\ Recalling that functions $d\omega^{0}, d\omega^{1}$ and $ d_{Y}\omega^{1}$ belongs to $\mathcal{C}_{c}^{\infty}(\overline{U}, \Lambda^{1}\mathcal{C}^{\infty}(Y))$, the function \\ $f(\cdot, \, ^{\sharp}d\omega^{0}+ \, ^{\sharp}d_{Y}\omega^{1}) $ belongs to $\mathcal{C}_{c}^{\infty}(\overline{U}, L^{\infty}(Y))$, so that (\ref{lem17}) implies 
	\begin{equation}\label{lem21}
	\lim_{\epsilon \to 0} \int_{U} f^{\epsilon}(\cdot,\, ^{\sharp}d\omega^{0}+ \, ^{\sharp}(d_{Y}\omega^{1})^{\epsilon}) vol_{p} = \iint_{U\times Y} f(p, \, ^{\sharp}d\omega^{0}_{p} + \, ^{\sharp}d_{Y}\omega^{1}_{(p,q)}) vol_{p}vol_{q}.
	\end{equation}
	On the other hand, since $d\omega^{\epsilon}_{p} = d\omega^{0}_{p} + \epsilon d\omega^{1}_{(p,p_{\epsilon})} + d_{Y}\omega^{1}_{(p,p_{\epsilon})}$, it follows from (\ref{lem19})
	\begin{equation}
	|f^{\epsilon}(\cdot, \, ^{\sharp}d\omega^{\epsilon}) - f^{\epsilon}(\cdot, \, ^{\sharp}d\omega^{0} + \,  ^{\sharp}(d_{Y}\omega^{1})^{\epsilon}) | \leq c\epsilon \quad \textup{in} \; U, \; \epsilon >0,
	\end{equation} 
	where $c = c\left( \|d\omega^{1}\|_{\infty}, \phi(\|d\omega^{0}\|_{\infty}), \phi(\|d_{Y}\omega^{1}\|_{\infty})  \right) >0$. Hence 
	\begin{equation}\label{lem20}
	f^{\epsilon}(\cdot, \, ^{\sharp}d\omega^{\epsilon}) - f^{\epsilon}(\cdot, \, ^{\sharp}d\omega^{0} + \, ^{\sharp}(d_{Y}\omega^{1})^{\epsilon}) \rightarrow 0 \quad \textup{in} \; L^{1}(U) \; \textup{as} \; \epsilon \to 0.
	\end{equation}
	Therefore we end the proof by combining (\ref{lem21})-(\ref{lem20}) with the decomposition 
	\begin{equation}
	\begin{array}{l}
	\int_{U}  f^{\epsilon}(\cdot, \, ^{\sharp}d\omega^{\epsilon})vol_{p} - \iint_{U\times Y} f(\cdot, \, ^{\sharp}d\omega^{0}+ \, ^{\sharp}d_{Y}\omega^{1}) vol_{p}vol_{q} \\ 
	=  \int_{U} \left[  f^{\epsilon}(\cdot, \, ^{\sharp}d\omega^{\epsilon}) -  f^{\epsilon}(\cdot, \, ^{\sharp}d\omega^{0}+ \, ^{\sharp}(d_{Y}\omega^{1})^{\epsilon}) 
	\right] vol_{p} + \int_{U} f^{\epsilon}(\cdot, \, ^{\sharp}d\omega^{0}+  \, ^{\sharp}(d_{Y}\omega^{1})^{\epsilon}) vol_{p} \\
	\;\;\;\; - \iint_{U\times Y} f(\cdot, \, ^{\sharp}d\omega^{0}+ \, ^{\sharp}d_{Y}\omega^{1}) vol_{p}vol_{q}. 
	\end{array}
	\end{equation}
\end{preuv}

We now define the new space for the homogenized integral functional. For this, let 
\begin{equation}
\mathbb{F}_{0}^{1}L^{\Phi} = H^{\Phi,d}_{0}(U)\times L^{\Phi}(U, \Lambda^{1}\, H^{\Phi, d}_{\#}(Y)). 
\end{equation}
We equip $\mathbb{F}_{0}^{1}L^{\Phi}$ with the norm 
\begin{equation}
\parallel \bm{\omega} \parallel_{\mathbb{F}_{0}^{1}L^{\Phi}} = \parallel d\omega^{0} \parallel_{L^{\Phi}(U,\Lambda^{1})} + \parallel d_{Y}\omega^{1} \parallel_{L^{\Phi}(U\times Y,\Lambda^{1,1})},\;\; \textup{with} \;\, \bm{\omega} = (\omega^{0}, \omega^{1}) \in \mathbb{F}_{0}^{1}L^{\Phi}.
\end{equation} 
With this norm, $\mathbb{F}_{0}^{1}L^{\Phi}$ is a Banach space admitting $F^{\infty}_{0} =   \mathcal{C}_{c}^{\infty}(U) \times \left[ \mathcal{C}_{c}^{\infty}(U)\otimes \mathcal{C}_{c}^{\infty}(U, \Lambda^{1}\mathcal{C}^{\infty}(Y)) \right]$ as a dense subspace. \\
Let now $\bm{\alpha} = (\alpha^{0}, \alpha^{1}) \in \mathbb{F}_{0}^{1}L^{\Phi}$, set $\mathbf{d}\bm{\alpha} = d\alpha^{0} + d_{Y}\alpha^{1} \in L^{\Phi}(U\times Y, \Lambda^{1,1})$ and define the functional $F$ on $\mathbb{F}_{0}^{1}L^{\Phi}$ by 
\begin{equation}
F(\bm{\alpha}) = \iint_{U\times Y} f(\cdot, \, ^{\sharp}\mathbf{d}\bm{\alpha})\,vol_{p}vol_{q},
\end{equation}
where\, $^{\sharp}\mathbf{d}\bm{\alpha} = \, ^{\sharp}d\alpha^{0} + \,  ^{\sharp}d_{Y}\alpha^{1}$ and the function $f$ here is defined as in corollary \ref{lem22}. \\
The hypotheses $(H_{1})-(H_{3})$ drive to the following lemma.
\begin{lemma}
	There exists a unique $\mathbf{u} = (u^{0}, u^{1}) \in \mathbb{F}_{0}^{1}L^{\Phi}$ such that 
	\begin{equation}\label{lem28}
	F(\mathbf{u}) = \left\{  \min F(\bm{\alpha})\, : \, \bm{\alpha} \in \mathbb{F}_{0}^{1}L^{\Phi} \right\}.
	\end{equation}
\end{lemma}


\subsection{Regularization} 

As in \cite{tacha1}, we regularize the integrand $f$ in order to get an approximating family of integrands $(f_{m})_{m\in\mathbb{N}^{\ast}}$ having in particular some properties $(H_{1})-(H_{3})$. Precisely, let $\chi_{m} \in \mathcal{C}^{\infty}_{0}(\mathbb{R}^{n})$ with $0\leq \chi_{m}$, $\textup{supp}\chi_{m} \subset \frac{1}{m}\overline{B_{n}}$ (where $\overline{B_{n}}$ denotes the closure of the open unit ball $B_{n}$ in $\mathbb{R}^{n}$) and $\int \chi_{m}(\eta)d\eta = 1$. Let $p\in U\subset M$ and $(p, \zeta_{p}) \in TU \cong U\times \mathbb{R}^{n}$, Setting 
\begin{equation}
f_{m}(p, \zeta_{p}) = \int \chi_{m}(\eta)f(p, \zeta_{p}-\eta)d\eta.
\end{equation}
The main properties of this new integrand are the following :
\begin{itemize}
	\item[$(H_{1})_{m}$] $f_{m}$ is measurable ;
	\item[$(H_{2})_{m}$] $f_{m}(p,\cdot)$ is strictly convex for almost all $p \in U$ ;
	\item[$(H_{3})_{m}$] There is constant $c_{5} > 0$ such that 
	\begin{equation}
	f_{m}(p,\zeta_{p}) \leq c_{5}(1+ \Phi(|\zeta_{p}|))
	\end{equation}
	 for almost all $p \in U$ and for all $\zeta_{p} \in T_{p}U$ ;
	\item[$(H_{4})_{m}$] $\dfrac{\partial f_{m}}{\partial \zeta_{p}}(p,\zeta_{p})$ exists for almost all $p \in U$ and for all $\zeta_{p} \in T_{p}U$, and there exists a constant $c_{6}=c_{6}(n) > 0$ such that 
	\begin{equation}
	\left| \dfrac{\partial f_{m}}{\partial \zeta_{p}}(p, \zeta_{p}) \right| \leq c_{6} (1+\Phi(|\zeta_{p}|)).
	\end{equation}
\end{itemize}
This being so, we obtain the results in proposition \ref{lem23} and in corollary \ref{lem24}, where $f$ is replaced by $f_{m}$. Moreover for every  $\mathbf{v} \in L^{\Phi}(U, \mathcal{C}^{\infty}(Y))^{n}$, as $m \to \infty$, one has
	\begin{equation}
	f_{m}(\cdot, \mathbf{v}) \rightarrow f(\cdot, \mathbf{v}) \quad \textup{in} \; L^{1}(U\times Y).
	\end{equation}
\noindent We are now ready to prove one of the most important results of this section.
\begin{prop}
	Let $(\mathbf{v}_{\epsilon})$ be a sequence in $ L^{\Phi}(U)^{n}$ which weakly 2s-converges (in each component) to $\mathbf{v} \in L^{\Phi}(U\times Y)^{n}$. Then for any integer $m\geq 1$, we have that there exists a constant $C'$ such that
	\begin{equation}
	\iint_{U\times Y} f_{m}(\cdot, \mathbf{v}) vol_{p}vol_{q} - \dfrac{C'}{m} \leq \lim_{\epsilon \to 0} \int_{U}  f(p^{\epsilon}, \mathbf{v}_{\epsilon}(p))  vol_{p}.
	\end{equation}
\end{prop}
\begin{preuv}
	Let integer $m\geq 1$ and let $(\mathbf{v}_{l})_{l}$ be a sequence in $ \left(\mathcal{C}^{\infty}_{c}(U)\otimes\mathcal{C}^{\infty}(Y)\right)^{n}$ such that $\mathbf{v}_{l} \rightarrow \mathbf{v}$ in  $L^{\Phi}(U\times Y)^{n}$ as $l\to \infty$. \\
	The convexity and the differentiability of $f_{m}(p,\cdot)$ imply (for any integer $l\geq 1$) 
	\begin{equation}
	\begin{array}{rcl}
	\int_{U}  f_{m}(p^{\epsilon}, \mathbf{v}(p,p^{\epsilon}))  vol_{p} & \geq & \int_{U}  f_{m}(p^{\epsilon}, \mathbf{v}_{l}(p,p^{\epsilon}))  vol_{p}  \\
	& &   + \int_{U} \dfrac{\partial f_{m}}{\partial \lambda} (p^{\epsilon}, \mathbf{v}_{l}(p,p^{\epsilon}))\cdot\left( \mathbf{v}(p,p^{\epsilon})  - \mathbf{v}_{l}(p,p^{\epsilon}) \right)  vol_{p} 
	\end{array}
	\end{equation}
	On the other hand, taking into account $(H_{1})_{m}$, $(H_{2})_{m}$ and $(H_{4})_{m}$, the function $p\rightarrow \frac{\partial f_{m}}{\partial \zeta_{p}}(\cdot, \mathbf{v}_{l}(p, \cdot))$ of $\overline{U}$ into $L^{\infty}(Y)^{n}$, denoted $\frac{\partial f_{m}}{\partial \zeta_{p}}(\cdot, \mathbf{v}_{l})$, belongs to $\mathcal{C}^{\infty}_{c}(\overline{U}, L^{\infty}(Y)^{n})$. Hence, arguing as in the first part of the proof of proposition \ref{lem23}, remark\ref{lem16} implies that 
	\begin{equation}
	\begin{array}{l}
	\lim_{\epsilon \to 0}	\iint_{U} \dfrac{\partial f_{m}}{\partial \zeta_{p}} (p^{\epsilon}, \mathbf{v}_{l}(p,p^{\epsilon}))\cdot \left( \mathbf{v}(p,p^{\epsilon}) - \mathbf{v}_{l}(p,p^{\epsilon}) \right)  dxd\mu  \\
	= \iint_{U\times Y} \dfrac{\partial f_{m}}{\partial \zeta_{p}} (q, \mathbf{v}_{l}(p,q))\cdot \left( \mathbf{v}(p,q) - \mathbf{v}_{l}(p,q) \right) vol_{p}vol_{q}. 
	\end{array} 
	\end{equation}
	Next, we observe that for a.e. $p \in \overline{U}$ and every $\zeta_{p} \in T_{p}U$ and a suitable positive constant $c'_{0}$ one has 
	\begin{equation}\label{bel1}
	f_{m}(p,\zeta_{p}) \leq f(p,\zeta_{p}) + \dfrac{1}{m}c'_{0} (1 + \phi(2(1+|\zeta_{p}|))).
	\end{equation}	
	Indeed, for a.e. $\omega$, every $\zeta_{p}$, $\nu_{p} \in T_{p}U$, by (\ref{lem7}), 
	\begin{equation}
	\begin{array}{lcr}
	f(p,\zeta_{p}) &  \leq &  f(p,\nu_{p}) +  c'_{0} \dfrac{ \Phi(2(1+ |\zeta_{p}| + |\nu_{p}|))}{1+ |\zeta_{p}| + |\nu_{p}|} |\zeta_{p} - \nu_{p}| \\
	& \leq & f(p,\nu_{p}) + c'_{0} (1 + \phi( 1+ |\zeta_{p}| + |\nu_{p}|))|\zeta_{p} - \nu_{p}|.
	\end{array}
	\end{equation}	
	Replacing $\zeta_{p}$ by $\zeta_{p}-\eta$ and $\nu_{p}$ by $\zeta_{p}$ respectively, we obtain : 
	\begin{equation}
	\begin{array}{lcr}
	f(p,\zeta_{p}-\eta) &  \leq &  f(p,\zeta_{p}) +  c'_{0} (1 + \phi( 1+ |\zeta_{p}-\eta| + |\zeta_{p}|))|\nu_{p}| \\
	& \leq & f(p,\zeta_{p}) + c'_{0} (1 + \phi( 1+ |\eta| + 2|\zeta_{p}|))|\nu_{p}|.
	\end{array}
	\end{equation}
	Let $m>0$, and assume $|\eta| \leq \frac{1}{m} \leq 1$, hence, 
	\begin{equation}
	f(p,\zeta_{p}-\eta)   \leq   f(p,\zeta_{p}) +  c'_{0} (1 + \phi(2( 1 + |\zeta_{p}|)))\frac{1}{m}. 
	\end{equation}	
	Multiplying both sides of the inequality, by $\chi_{m}$, we get : 
	\begin{equation}
	f(p,\zeta_{p}-\eta)\chi_{n}(\eta)   \leq   f(p,\zeta_{p})\chi_{n}(\eta) +  \frac{1}{m}c'_{0} (1 + \phi(2( 1 + |\zeta_{p}|))) \chi_{n}(\eta). 
	\end{equation}
	Integration leads to (\ref{bel1}). Hence, given $\mathbf{v}_{\epsilon}(p)= \mathbf{v}(p, p^{\epsilon})$, we have 	
	\begin{equation}
	f_{m}(p^{\epsilon},\mathbf{v}_{\epsilon}) \leq f(p^{\epsilon},\mathbf{v}_{\epsilon}) + \dfrac{1}{m} c'_{0} (1 + \phi(2(1+|\mathbf{v}_{\epsilon}|))).
	\end{equation}
	Thus 
	\begin{equation}
	\begin{array}{rcl}
	\int_{U} f_{m}(p^{\epsilon},\mathbf{v}_{\epsilon}) vol_{p} & \leq & \int_{U} f(p^{\epsilon},\mathbf{v}_{\epsilon}) vol_{p}  + \dfrac{1}{m}C|U|  \\
	& & + \dfrac{c'_{0}}{m} \int_{U} \alpha \dfrac{\phi(2(1+|\mathbf{v}_{\epsilon}|))}{\alpha} vol_{p},
	\end{array}
	\end{equation}
	with $0 < \alpha \leq 1$.
	But 
	\begin{equation}
	\alpha \dfrac{\phi(2(1+|\mathbf{v}_{\epsilon}|))}{\alpha} \leq \widetilde{\Phi}(\alpha\phi(2(1+|\mathbf{v}_{\epsilon}|))) + \Phi\left(\frac{1}{\alpha}\right) \leq \alpha\widetilde{\Phi}(\phi(2(1+|\mathbf{v}_{\epsilon}|))) + \Phi\left(\frac{1}{\alpha}\right).
	\end{equation}
	Set $\Gamma_{1} = \{ (p \in U : 2(1+|\mathbf{v}_{\epsilon}|) > t_{0} \}$, $\Gamma_{2} = U \backslash \Gamma_{1}$.  \\
	Hence,we get 
	\begin{equation}
	\begin{array}{rcl}
	\int_{U} \alpha \dfrac{\phi(2(1+|\mathbf{v}_{\epsilon}|))}{\alpha} vol_{p} & \leq & \int_{U} \alpha\widetilde{\Phi}(\phi(2(1+|\mathbf{v}_{\epsilon}|))) vol_{p} + \Phi\left(\frac{1}{\alpha}\right) |U| \\
	& \leq &  |\Gamma_{2}|\alpha \widetilde{\Phi}(\phi(t_{0})) +  \Phi\left(\frac{1}{\alpha}\right) |U| + \alpha \int_{\Gamma_{1}} \Phi(4(1+|\mathbf{v}_{\epsilon}|)) vol_{p}.
	\end{array}
	\end{equation}
	Let $C > 1 + \|4(1+|\mathbf{v}_{\epsilon}|)\|_{\Phi,U}$. Then $\int_{U} \Phi\left(\frac{4(1+|\mathbf{v}_{\epsilon}|)}{C}\right) vol_{p} \leq 1$.  \\
	Since $$\Phi(4(1+|\mathbf{v}_{\epsilon}|)) = \Phi\left(C\frac{4(1+|\mathbf{v}_{\epsilon}|)}{C}\right) \leq K(C) \Phi\left(\frac{4(1+|\mathbf{v}_{\epsilon}|)}{C}\right) \; \textup{whenever}\; \frac{4(1+|\mathbf{v}_{\epsilon}|)}{C} \geq t_{0}$$. \\
	Set $\Gamma_{3} = \left\{ p \in \Gamma_{1} : \frac{4(1+|\mathbf{v}_{\epsilon}|)}{C} \geq t_{0} \right\}$, $\Gamma_{4} = \Gamma_{1}\backslash\Gamma_{3}$. \\
	Hence,
	\begin{equation}
	\begin{array}{rcl}
	\int_{\Gamma_{1}} \Phi(4(1+|\mathbf{v}_{\epsilon}|)) vol_{p} & = & \int_{\Gamma_{4}} \Phi(4(1+|\mathbf{v}_{\epsilon}|)) vol_{p}  + \int_{\Gamma_{3}} \Phi(4(1+|\mathbf{v}_{\epsilon}|)) vol_{p} \\
	& \leq &  |\Gamma_{4}|\Phi(C t_{0}) + \int_{\Gamma_{3}} \Phi(4(1+|\mathbf{v}_{\epsilon}|)) vol_{p}  \\
	& \leq &  |\Gamma_{4}|\Phi(C t_{0}) + \alpha \int_{\Gamma_{3}} \Phi\left(C\dfrac{4(1+|\mathbf{v}_{\epsilon}|)}{C}\right) vol_{p} \\
	& \leq & |\Gamma_{4}|\Phi(C t_{0}) + K(C) \int_{\Gamma_{3}} \Phi\left(\frac{4(1+|\mathbf{v}_{\epsilon}|)}{C}\right) vol_{p}  \\
	& \leq & |\Gamma_{4}|\Phi(C t_{0}) + K(C) \int_{U} \Phi\left(\frac{4(1+|\mathbf{v}_{\epsilon}|)}{C}\right) vol_{p}.
	\end{array}
	\end{equation}
	Since $\Phi \in \Delta_{2}$, and $(\mathbf{v}_{\epsilon})$ is bounded in $L^{\Phi}(U)^{n}$ it results that $\int_{U} \Phi\left(4(1+|\mathbf{v}_{\epsilon}|)\right) vol_{p}$ is also bounded.  
	Then we have 
	\begin{equation}
	\begin{array}{rcl}
	\int_{U} f_{m}(p^{\epsilon},\mathbf{v}_{\epsilon}) vol_{p} & \leq & \int_{U} f(p^{\epsilon},\mathbf{v}_{\epsilon}) vol_{p} + \dfrac{1}{n}C|U|  \\
	& &  + \dfrac{c'_{0}}{m} \bigg( \alpha |U|\widetilde{\Phi}(\phi(t_{0})) + \Phi\left(\frac{1}{\alpha}\right)|U|  \\
	& & + \alpha\left( |\Gamma_{4}|\Phi(C t_{0}) + K(C) \right) \int_{U} \Phi \left(  \dfrac{4(1+|\mathbf{v}_{\epsilon}|)}{C} \right) vol_{p}    \bigg)  \\
	& \leq & \int_{U} f(p^{\epsilon},\mathbf{v}_{\epsilon}) vol_{p} + \dfrac{1}{m}C'  
	\end{array}
	\end{equation}
	for a suitably big constant $C'$. \\
	Thus
	\begin{equation}
	\begin{array}{l}
	\liminf_{\epsilon \to 0} \int_{U}  f(p^{\epsilon}, \mathbf{v}(p,p^{\epsilon}))  vol_{p}  \\
	\geq  \iint_{U\times Y}  f_{m}(p, \mathbf{v}_{l}(p,q))  vol_{p}vol_{q} \\
	\quad - \frac{C'}{m}  + \iint_{U\times Y} \dfrac{\partial f_{m}}{\partial \zeta_{p}} (q, \mathbf{v}_{l}(p,q))\cdot \left( \mathbf{v}(p,q) - \mathbf{v}_{l}(p,q) \right)  vol_{p}vol_{q}. 
	\end{array}
	\end{equation}	
	Using $(H_{4})_{m}$ combined with the H\"{o}lder's inequality and (\ref{lem11}) yields 
	\begin{equation}
	\left| \iint_{U\times Y} \dfrac{\partial f_{m}}{\partial \zeta_{p}} (q, \mathbf{v}_{l}(p,q))\cdot \left( \mathbf{v}(p,q) - \mathbf{v}_{l}(p,q) \right)  vol_{p}vol_{q} \right| \leq c'_{0} \, \|1 + \phi(|\mathbf{v}_{l}|)\|_{\widetilde{\Phi},U\times Y}\cdot \|\mathbf{v}-\mathbf{v}_{l}\|_{\Phi,U\times Y}.
	\end{equation}
	Since $\mathbf{v}_{l} \rightarrow \mathbf{v}$ in  $L^{\Phi}(U\times Y)^{n}$ as $l\to \infty$, it follows that for $\delta > 0$ arbitrarily fixed, there exists $l_{0} \in \mathbb{N}$ such that 
	\begin{equation}
	\left| \iint_{U\times Y} \dfrac{\partial f_{m}}{\partial \zeta_{p}} (p, \mathbf{v}_{l}(p,q))\cdot \left( \mathbf{v}(p,q) - \mathbf{v}_{l}(p,q) \right)  vol_{p}vol_{q} \right| \leq \delta 
	\end{equation}
	for all $l \geq l_{0}$.
	Hence for all $l \geq l_{0}$,
	\begin{equation}
	\lim_{\epsilon \to 0} \int_{U}  f(p^{\epsilon}, \mathbf{v}_{\epsilon})  vol_{p} \geq \iint_{U\times Y} f_{m}(p, \mathbf{v}_{l}(p,q)) vol_{p}vol_{q} - \delta - \dfrac{C'}{m}.
	\end{equation}
	Now sending $l \to \infty$ we have 
	\begin{equation}
	\lim_{\epsilon \to 0} \int_{U}  f(p^{\epsilon}, \mathbf{v}_{\epsilon})  vol_{p} \geq \iint_{U\times Y} f_{m}(q, \mathbf{v}(p,q))  vol_{p}vol_{q} - \delta - \dfrac{C'}{m}.
	\end{equation}
	Since $\delta$ is arbitrarily fixed, we are led at once to the result by letting $\delta \to 0$.
\end{preuv}

Letting $m\to\infty$, and replacing $\mathbf{v}_{\epsilon}$ by\,  $^{\sharp}du^{\epsilon}$, with $du^{\epsilon}$ weakly two-scale converges  to $du = du_{0} + d_{Y}u_{1}$ in $ L^{\Phi}(U\times Y, \Lambda^{1,0})$, one obtains the following result : 


\begin{cor}\label{lem32}
	Let $(u^{\epsilon})$ be a sequence in $H^{\Phi,d}(U)$. Assume that $(du^{\epsilon})$  weakly two-scale converges to $\textup{\textbf{d}}\textup{\textbf{u}} = du^{0} + d_{Y}u^{1}$ in  $  L^{\Phi}(U\times Y, \Lambda^{1,0})$, where $\mathbf{u} = (u^{0}, u^{1}) \in \mathbb{F}^{1}_{0}L^{\Phi}$. We denotes \, $^{\sharp}\textup{\textbf{d}}\textup{\textbf{u}} = \, ^{\sharp}du^{0} + \, ^{\sharp}d_{Y}u^{1}$ and we have  
	\begin{equation}
	\iint_{U\times Y} f(p, \,  ^{\sharp}\mathbf{d}\mathbf{u}_{(p,q)}) vol_{p}vol_{q} \leq \lim_{\epsilon \to 0} \int_{U}  f(p^{\epsilon}, \, ^{\sharp}du^{\epsilon}_{p})  vol_{p}
	\end{equation}
\end{cor}

\subsection{Main homogenization result}

Our main objective in this subsection is to prove the following 
\begin{theorem}\label{lem37}
	For each $\epsilon > 0$, let $(u^{\epsilon}) \in H^{\Phi,d}_{0}(U)$ be the unique solution of (\ref{lem26}). Then, as $\epsilon \to 0$, 
	\begin{equation}\label{lem30}
	u^{\epsilon} \;\;\; \substack{2s \\ \rightharpoonup} \;\;\; u^{0} \quad \textup{in} \;\, L^{\Phi}(U)  
	\end{equation}
	and 
	\begin{equation}\label{lem31}
	du^{\epsilon} \;\;\; \substack{2s \\ \rightharpoonup} \;\;\; \textup{\textbf{d}}\textup{\textbf{u}}= du^{0} + d_{Y}u^{1} \quad \textup{in} \;\, L^{\Phi}(U\times Y, \Lambda^{1,0}),  
	\end{equation}
	where $\mathbf{u} = (u^{0}, u^{1}) \in \mathbb{F}_{0}^{1}L^{\Phi}$ is the unique solution to the minimization problem (\ref{lem28}). 
\end{theorem}
\begin{preuv}
	In view of the growth conditions in $(H_{3})$, the sequence $(u^{\epsilon})_{\epsilon\in E}$ is bounded in $H^{\Phi,d}_{0}(U)$ and so the sequence $(f^{\epsilon}(\cdot, \,  ^{\sharp}du^{\epsilon}))_{\epsilon>0}$ is bounded in $L^{1}(U)$. Thus, given an arbitrary fundamental sequence $E$, we get by theorem \ref{math5} the existence of a subsequence $E'$ from $E$ and a couple $\mathbf{u} = (u^{0}, u^{1}) \in \mathbb{F}_{0}^{1}L^{\Phi}$ such that (\ref{lem30})-(\ref{lem31}) hold when $E' \ni \epsilon \to 0$. The sequence $(F_{\epsilon}(u^{\epsilon}))_{\epsilon>0}$ consisting of real numbers being bounded, since $(u^{\epsilon})_{\epsilon>0}$ is bounded in $H^{\Phi,d}_{0}(U)$, there exists a subsequence from $E'$ not relabeled such that $\lim_{E' \ni\epsilon \to 0} F_{\epsilon}(u^{\epsilon})$ exists. \\
	It remains to verify that $\mathbf{u} = (u^{0}, u^{1})$ solves (\ref{lem28}). In fact, if $\mathbf{u}$ solves this problem, then thanks to the uniqueness of the solution of (\ref{lem28}), the whole sequence $(u^{\epsilon})_{\epsilon>0}$ will verify (\ref{lem30}) and (\ref{lem31}) when $\epsilon \to 0$. Thus our only concern here is to check that $\mathbf{u}$ solves problem (\ref{lem28}). To this end, in view of corollary \ref{lem32}, we have 
	\begin{equation}\label{lem34}
	\iint_{U\times Y} f(p, \,  ^{\sharp}\mathbf{d}\mathbf{u}) vol_{p}vol_{q} \leq \lim_{E'\ni\epsilon \to 0} \int_{U}  f(p^{\epsilon}, du^{\epsilon}_{p})  vol_{p}.
	\end{equation}
	On the other hand, let us establish an upper bound for 
	\begin{equation}
	\int_{U} f(p^{\epsilon}, \, ^{\sharp}du^{\epsilon}_{p})  vol_{p}.
	\end{equation}
	To do that, let $\phi = \left(\psi^{0}, \psi^{1}\right) \in F^{\infty}_{0}$ with $\psi_{0} \in \mathcal{C}_{c}^{\infty}(U)$, $\psi_{1} \in \mathcal{C}_{c}^{\infty}(U, \Lambda^{1})\otimes \mathcal{C}^{\infty}(Y)$. 	Define $\phi^{\epsilon}$ as in corollary \ref{lem24}. Since $u^{\epsilon}$ is the minimizer, one has 
	\begin{equation}
	\int_{U} f(p^{\epsilon}, \, ^{\sharp}du^{\epsilon}_{p})  vol_{p}  \leq \int_{U}  f(p^{\epsilon}, \, \phi^{\epsilon}_{(p,p^{\epsilon})})  vol_{p}.
	\end{equation}
	Thus, using corollary \ref{lem24} we get 
	\begin{equation}
	\lim_{E'\ni\epsilon \to 0} \int_{U} f(p^{\epsilon}, \, ^{\sharp}du^{\epsilon}_{p})  vol_{p} \leq\iint_{U\times Y}  f(\cdot, \, ^{\sharp}d\psi^{0} + \, ^{\sharp}d_{Y}\psi^{1})  vol_{p}vol_{q},
	\end{equation}
	for any $\phi \in F^{\infty}_{0}$, and by density, for all $\phi \in \mathbb{F}_{0}^{1}L^{\Phi}$. From which we get 
	\begin{equation}\label{lem33}
	\lim_{E'\ni\epsilon \to 0} \int_{U} f(p^{\epsilon}, \, ^{\sharp}du^{\epsilon}_{p})  vol_{p} \leq  \inf_{\bm{\alpha} \in \mathbb{F}_{0}^{1}L^{\Phi}} \iint_{U\times Y}  f(\cdot, \, ^{\sharp}\mathbf{d}\bm{\alpha})  vol_{p}vol_{q}.
	\end{equation}
	Inequalities (\ref{lem34}) and (\ref{lem33}) yield 
	\begin{equation}
	\iint_{U\times Y} f(p, \, ^{\sharp}\mathbf{d}\mathbf{u}) vol_{p}vol_{q} = \inf_{\bm{\alpha} \in \mathbb{F}_{0}^{1}L^{\Phi}} \iint_{U\times Y}  f(\cdot, \, ^{\sharp}\mathbf{d}\bm{\alpha})  vol_{p}vol_{q}
	\end{equation}
	i.e. (\ref{lem28}). The proof is complete.
\end{preuv}

\begin{remark}
	Indeed when we consider the particular case where the Riemannian manifolds $M = \mathbb{R}^{n}$, $Y= [0,1]^{n}$ (the unit cube of $\mathbb{R}^{n}$), $U = \Omega$ an open bounded subset of $\mathbb{R}^{n}$ then one can write a point $p^{\epsilon} \in \Omega$ as $p^{\epsilon} = x^{\epsilon} = \frac{x}{\epsilon}$ where $x\in \mathbb{R}^{n}$. Moreover we can write any function $f : TU\rightarrow \mathbb{R}$ as $f : \mathbb{R}^{n}\times \mathbb{R}^{n}\rightarrow \mathbb{R}$ and we identified for any function $\omega : \Omega\rightarrow \mathbb{R}$,  $^{\sharp}d\omega(x) \cong \nabla\omega(x)$, $H^{\Phi,d}(\Omega) \cong W^{1}L^{\Phi}(\Omega) = \left\{ u\in L^{\Phi}(\Omega) : \dfrac{\partial u}{\partial x_{i}} \in L^{\Phi}(\Omega), \; 1\leq i \leq n \right\}$. 
	 Therefore, the homogenization problem (\ref{lem26}) is equivalent to the periodic problem (see \cite{tacha1}),
	\begin{equation}
	\min \left\{ \int_{\Omega} f\left(\frac{x}{\epsilon}, \nabla \omega(x) \right)dx \; : \; \omega \in  W^{1}_{0}L^{\Phi}(\Omega) \right\}.
	\end{equation}
\end{remark}

\vspace{0.5cm}

 

\vspace{1cm} 

\par \textit{Current adress :} $^{\ddagger }$The University of Maroua, Department of Mathematics and Computer Science, P.O. Box 814 Maroua, Cameroon.
\par \textit{E-mail :} takougoumfranckarnold@gmail.com 
 
\vspace{0.3cm}

\textit{Current adress :} $^{\dagger }$The University of Bamenda, Higher Teachers Training
College, Department of Mathematics, P.O. Box 39, Bambili Cameroon 
\par \textit{E-mail :} adressfotsotachago@yahoo.fr 
 
\vspace{0.3cm}

\par\textit{Current adress :} $^{\dagger ^{\dagger }}$The University of Maroua, Department of Mathematics and Computer Science, P.O. Box 814 Maroua, Cameroon. 
\par \textit{E-mail :} joseph.dongho@fs.univ-maroua.cm   


\end{document}